\documentclass[pdflatex,sn-mathphys,Numbered]{sn-jnl}

\usepackage{hyperref}
\usepackage{subcaption}
\usepackage{fullpage}
\usepackage{graphicx}%
\usepackage{multirow}%
\usepackage{amsmath,amssymb,amsfonts}%
\usepackage{amsthm}%
\usepackage{mathrsfs}%
\usepackage[title]{appendix}%
\usepackage{xcolor}%
\usepackage{textcomp}%
\usepackage{manyfoot}%
\usepackage{booktabs}%
\usepackage{algorithm}%
\usepackage{algorithmicx}%
\usepackage{algpseudocode}%
\usepackage{listings}%
\theoremstyle{thmstyleone}%
\newtheorem{theorem}{Theorem}
\newtheorem{lemma}[theorem]{\bf Lemma}

\theoremstyle{thmstyletwo}%
\newtheorem{remark}{Remark}%

\theoremstyle{thmstylethree}%
\newtheorem{definition}{Definition}%
\newtheorem{assumption}{Assumption}

\DeclareMathOperator{\rank}{rank}

\newcommand{\R}{\mathbb{R}}

\newcommand{\Rnn}{\mathbb{R}^{n \times n}}
\newcommand{\Rmn}{\mathbb{R}^{m \times n}}

\newcommand{\Rnr}{\mathbb{R}^{n \times r}}

\newcommand{\Rmr}{\mathbb{R}^{m \times r}}

\newcommand{\Rrr}{\mathbb{R}^{r \times r}}

\DeclareMathOperator{\myspan}{span}
\newcommand{\norm}[1]{\left\| #1 \right\|}
\newcommand{\scal}[1]{\left\langle #1 \right\rangle}

\newcommand{\Mr}{\mathcal{M}_r}

\newcommand{\Tr}{\mathcal{T}_r}

\newcommand{\proj}[2]{\mathcal{P}_{#1} \left[ #2 \right]}

\newcommand{\Ret}[1]{\mathcal T_r \left( #1 \right) }

\newcommand{\vphi}[1]{\varphi_1 \left( #1 \right)}
\newcommand{\cL}{\mathcal L}
\newcommand{\cG}{\mathcal G}
\newcommand{\cF}{\mathcal F}

\newcommand{\ellF}{\ell_{\mathcal F}}
\newcommand{\Yrungehalf}{Y_{1/2}^{\mathrm{PR}}}

\usepackage{accents}
\newcommand*{\dt}[1]{%
  \accentset{\mbox{\large\bfseries .}}{#1}}

\begin{document}

\title[Projected exponential methods]{Projected exponential methods for stiff dynamical low-rank approximation problems}


\author*[1]{\fnm{Benjamin} \sur{Carrel}}\email{benjamin.carrel@unige.ch}

\author[1]{\fnm{Bart} \sur{Vandereycken}}\email{bart.vandereycken@unige.ch}


\affil*[1]{\orgdiv{Mathematics}, \orgname{University of Geneva}, \orgaddress{\street{Rue du Conseil-Général}, \city{Geneva}, \postcode{1205}, \country{Switzerland}}}

\abstract{The numerical integration of stiff equations is a challenging problem that needs to be approached by specialized numerical methods. 
Exponential integrators form a popular class of such methods since they are provably robust to stiffness and have been successfully applied to a variety of problems. 
The dynamical low-rank approximation is a recent technique for solving high-dimensional differential equations by means of low-rank approximations.
However, the domain is lacking numerical methods for stiff equations since existing methods are either not robust-to-stiffness or have unreasonably large hidden constants.

In this paper, we focus on solving large-scale stiff matrix differential equations with a Sylvester-like structure, that admit good low-rank approximations.
We propose two new methods that have good convergence properties, small memory footprint and that are fast to compute.
The theoretical analysis shows that the new methods have order one and two, respectively.
We also propose a practical implementation based on Krylov techniques. 
The approximation error is analyzed, leading to a priori error bounds and, therefore, a mean for choosing the size of the Krylov space.
 Numerical experiments are performed on several examples, confirming the theory and showing good speedup in comparison to existing techniques.}



\maketitle

\section{Introduction} \label{sec:introduction}

In this work, we are interested in the numerical integration of large-scale and stiff initial value problems. A stiff problem is a problem for which an implicit method, like implicit Euler, is significantly faster than its explicit counterpart, like explicit Euler. Since stiffness is intrinsic to the problem, it cannot be avoided by reformulation and needs to be integrated by specialized numerical methods. The archetypal example is the heat equation $\partial_t u = \Delta u + g(u)$. 

We consider discretized partial differential equations (PDEs), where the stiff components have been approximated by linearization, leading to a Sylvester-like differential equation 
$$\dt{X}(t) = A X(t) + X(t) B + \mathcal G(t, X(t)), \quad X(0) = X_0 \in \Rmn.$$ 
The time independent matrices $A \in \R^{m \times m}$ and $B \in \Rnn$ are typically very large but also sparse or structured. When we are dealing with such large-scale problems, we cannot afford to store nor compute the dense solution $X(t)$. We therefore need some kind of data sparse representation for $X(t)$, like low-rank approximations. Applications of such large Sylvester-like differential equations arise in various fields, including numerical simulations, control theory, signal processing; see also Section~\ref{sec:numerical_experiments} and references therein for examples and more details. 

The dynamical low-rank approximation (DLRA) \cite{koch2007dynamical} is a popular technique that allows to compute an approximation of the solution by a low-rank matrix: $Y(t) \approx X(t)$ such that $\rank(Y(t)) = r \ll m, n$. The technique has been successfully applied in various areas, including control theory, signal processing, machine learning, image compression, and quantum physics. The projector-splitting integrator \cite{lubich2014projector} as well as the so-called unconventional integrator \cite{ceruti2022unconventional} are two efficient integrators of the DLRA, but the current convergence bounds \cite{kieri2016discretized} rely on the Lipschitz constant of the vector field, making these methods theoretically and practically not suitable for stiff problems.

The key idea of exponential methods \cite{hochbruck2010exponential} is to solve exactly the stiff part of the ODE. These methods are robust to stiffness by construction and their theoretical properties are well studied. In the same spirit, a low-rank Lie--Trotter splitting \cite{ostermann2019convergence} is shown to be robust to stiffness and can be used to solve the DLRA for parabolic problems. However, the Lie--Trotter and Strang splittings employed there suffer from large hidden constants in their error estimates, making the methods unfortunately unpractical. 

In this paper, we propose a new class of methods that we call projected exponential methods. These methods are robust to stiffness, have small hidden constants in the error estimates, and the computations preserve the low-rank structure of the solution. The growth in rank of the sub-steps is contained by means of projections and Krylov iterative methods. 

\section{Notations and preliminaries} \label{sec:preliminaries}

Let $A \in \R^{m \times m}$ and $B \in \Rnn$ be two matrices which lead to the stiffness\footnote{
In this work, we call a method robust to stiffness if it does not depend on the Lipschitz constant of the stiff part of the vector field. It may however depend on the one-sided Lipschitz constant.}. 
One can think of these matrices as discretized operators obtained from a particular PDE. Under mesh refinement, their Lipschitz constants will grow in a unbounded way, causing stiffness.  
In this work, we are interested in solving general Sylvester-like differential equations of the form
\begin{equation} \label{eq:full problem}
    \begin{aligned}
        \dt{X}(t) & = A X(t) + X(t) B + \cG(X(t)), \quad t \in [0, T] \\
        X(0)      & = X_0 \in \Rmn,
    \end{aligned}
\end{equation}
where $\cG \colon \Rmn \rightarrow \Rmn$ is a non-linear but
non-stiff\footnote{By non-stiff, we mean that the Lipschitz constant of $\cG$ remains bounded under mesh refinement. In practice, it should also not be too large.} operator. 
For ease of notation, we will often denote the linear Sylvester operator as $\mathcal L$, defined by $\mathcal LX = AX + XB$.
Throughout, we will denote $f(t) = \cG(X(t))$ and assume that $f(t)$ is absolutely continuous on the compact domain of integration $[0,T]$. 
It is sometimes assumed to be higher-order differentiable when needed. 
Without loss of generality, equation~\eqref{eq:full problem} is autonomous. 
Occasionally, we will denote the full vector field $\mathcal F = \mathcal L + \mathcal G$. Whenever we use calligraphic notation for an operator, it has to be interpreted as acting on the matrix: $\mathcal L X = \mathcal L(X)$ is not (in general) the matrix multiplication of $\mathcal L$ and $X$ but the linear operator $\mathcal L$ applied to the matrix $X$. The same holds for matrix functions, like $e^{t \cL}$ in the next lemma.
\begin{lemma}[Closed-form solution] \label{lemma:closed form solution}
    The solution to the initial value problem
    $$\dt{X}(t) = \cL X(t)  + \cG(X(t)), \quad X(0) = X_0,$$
    is given by the closed-form formula
    \begin{equation}\label{eq:closed form solution integral}
    X(t) = e^{t \cL} X_0 + \int_0^t e^{(t-s) \cL}  \cG(X(s))  ds.
    \end{equation}
\end{lemma}
The following definitions will be useful in the analysis of the methods. Throughout, we will denote by $\| \cdot \|$ the Frobenius norm, which is induced by the Euclidean inner product $\langle \,\cdot\,, \,\cdot\, \rangle$. Other norms will be specified explicitly.
\begin{definition}[One-sided Lipschitz] \label{def: one sided Lipschitz}
    Let $\cF\colon \mathbb{R}^{m \times n} \to \mathbb{R}^{m \times n}$ be an operator, not necessarily linear.
    Its one-sided Lipschitz constant $\ell_\cF$ is the smallest constant such that 
    \begin{equation}\label{eq:def_one-sided-Lip}
        \scal{X-Y, \cF(X) - \cF(Y)} \leq \ell_\cF \ \|X-Y\|^2 \qquad \forall X,Y \in \mathbb{R}^{m \times n}.
    \end{equation}
\end{definition}
For linear operators, like $\cL$, the one-sided Lipschitz constant can be calculated explicitly. In particular, condition~\eqref{eq:def_one-sided-Lip} is equivalent to
\[
 \ell_\cL = \max_{X-Y \neq 0} \frac{\scal{X-Y,\cL (X-Y)}}{\|X-Y\|^2} = \max_{W \neq 0} \frac{\scal{W,\cL W}}{\|W\|^2}.
\]
Since $\cL$ is defined\footnote{While we focus on $\R$, the following derivations are similar in $\mathbb{C}$.} 
on $\mathbb{R}$, we have $\scal{W,\cL W} = \scal{W,\cL^* W}$ with $\cL^*$ the adjoint of $\cL$. This allow us to write the one-sided Lipschitz constant as the maximal value of the Rayleigh quotient of a symmetric matrix, which is attained by the maximal eigenvalue of this matrix:
\[
\ell_\cL = \tfrac{1}{2} \max_{W \neq 0}  \frac{\scal{W,(\cL + \cL^*) W}}{\| W \|^2} = \tfrac{1}{2}\lambda_{\max}(\cL + \cL^*).
\]
The above formula is valid for any linear $\cL$. In case, $\mathcal L(X) = AX + XB$ is the Sylvester operator, we can go one step further and obtain
\[
 \ell_\cL = \tfrac{1}{2} \lambda_{\max}(A+A^*) + \tfrac{1}{2} \lambda_{\max}(B+B^*).
\]
This follows from the matrix representation of $\cL$ as $I \otimes A + B^T \otimes I$ and properties of the Kronecker product; see, e.g., \cite[Chap.~4.4]{hornTopicsMatrixAnalysis1991}.

\begin{definition}[Lipschitz continuity] \label{def: Lipschitz}
    An operator $\cF : \Rmn \rightarrow \Rmn$ is locally Lipschitz-continuous in a strip along the exact solution $X$ of \eqref{eq:full problem} if there exists $R > 0$ and $L_{\cF} \in \R$ such that, for all $t \in [0,T]$, $$\norm{\cF(Y(t)) - \cF(Z(t))} \leq L_\cF \norm{Y(t) - Z(t)},$$ if $\norm{X(t)-Y(t)} \leq R$ and $\norm{X(t)-Z(t)} \leq R$.
\end{definition}

\begin{remark}\label{remark: bound ellF sum}
    By elementary properties of the one-sided and the standard Lipschitz constants, we have for $\cF = \cL + \cG$ that
    \[
     \ell_\cF \leq \ell_\cL + \ell_\cG \leq \ell_\cL + L_\cG.
    \] 
    For diffusive problems, $\ell_\cL$ can be negative. As we will see,  obtaining error bounds in terms of one-sided Lipschitz constants can therefore be very beneficial in this case, even when $L_\cG$ is always positive.
\end{remark}

For completeness, we briefly introduce $\varphi$-functions which are central in the analysis and implementation of exponential integrators. 
\begin{definition}[$\varphi$-functions \cite{skaflestad2009scaling}] \label{def: phi functions} %
Let $\cL\colon \R^{m \times n} \to \R^{m \times n}$ be a linear operator. The $\varphi$-functions are defined as
\begin{align*}
&\varphi_0 (t \mathcal L) = e^{t \mathcal L}, \\
&\varphi_k (t \mathcal L) = \frac{1}{t^k} \int_0^t e^{(t-s) \mathcal L} \frac{s^{k-1}}{(k-1)!} ds, \quad k \geq 1.
\end{align*}
\end{definition}
A number of equivalent definitions exist for $\varphi$-functions. We first mention 
\begin{equation}\label{eq: phi truncated Taylor exp}
 \varphi_k(z) = \frac{e^z - p_{k-1}(z)}{z^k}, \qquad p_{k-1}(z) = \sum_{j=0}^{k-1} \frac{z^j}{j!},
\end{equation}
for scalar values $z\neq 0$ which is useful for explicit calculations with Taylor series. Next, the following definition makes the link with ODEs clear.
\begin{lemma}[$\varphi$-functions as IVPs \cite{skaflestad2009scaling,gockler2014uniform}] \label{lemma: phi functions as IVPs} %
Let $Z_i \in \Rmn$ for $i = 0, \ldots, n$. The $\varphi$-functions are equivalently defined by the differential equations
\begin{align*}
&Z(t) = \varphi_0 (t \mathcal L) Z_0 && \iff \quad \dt{Z}(t) = \mathcal L Z(t), \quad Z(0) = Z_0, \\
&Z(t) = t^k \varphi_k (t \mathcal L) Z_k && \iff \quad \dt{Z}(t) = \mathcal L Z(t) + \frac{t^{k-1}}{(k-1)!} Z_k, \quad Z(0) = 0.
\end{align*}
By linearity, we have that
\begin{equation*}
Z(t) = e^{t \mathcal L} Z_0 + \sum_{k=1}^n t^k \varphi_k(t \mathcal L) Z_k \iff \dt{Z}(t) = \mathcal L Z(t) + \sum_{k=1}^n \frac{t^{k-1}}{(k-1)!} Z_k, \quad Z(0) = Z_0.
\end{equation*}
\end{lemma}

Bounding these $\varphi$-functions appropriately is crucial for the analysis of exponential integrators introduced in the next section.
The derivation of such bounds requires to use the logarithmic norm (which, despite the name, is not a norm).

\begin{definition}[Logarithmic norm] \label{def: logarithmic norm}
    The logarithmic norm of a linear operator $\cL \colon \Rmn \rightarrow \Rmn$ is defined as 
    $$\mu_{\cL}~=~\lim_{h \rightarrow 0^+} \frac{\norm{I + h\cL}_2 - 1}{h},$$
    where $\| \, \cdot \, \|_2$ is the spectral norm.
\end{definition}
It is well known (see, e.g., \cite[Thm.~I.10.5]{hairersolving}) that the logarithmic norm for the spectral norm satisfies
\[
 \mu_{\cL} = \tfrac{1}{2}\lambda_{\max}(\cL+\cL^*).
\]
This implies $\mu_{\cL} = \ell_\cL$, the logarithmic norm equals the one-sided Lipschitz constant for a linear operator $\cL$. From now on, we will only use the notation $\ell = \ell_\cL$, and specify the index when needed.

\begin{lemma}[$\varphi$-functions are bounded] \label{lemma: phi functions are bounded}
    For any matrix $X \in \Rmn$, we have for all $k=0,1,\ldots$ that
    \begin{align*}
        &\norm{\varphi_k(h \cL) X}_{\star} \leq \varphi_k(h\ell_\cL) \norm{X}_{\star}, \quad h \geq 0,
    \end{align*}
    where $\norm{\, \cdot\,}_{\star}$ is either the spectral norm $\norm{\,\cdot\,}_2$ or the Frobenius norm $\norm{\,\cdot\,}$.
\end{lemma}
\begin{proof}
    Consider the differential equation $\dt{Z}(t) = \cL Z(t) + R(t)$, where $R(t) \in \mathcal C^2$ is specified below. Denote $\ell=\ell_\cL$. Then, by definition of the Dini derivative, we have
    $$\begin{aligned}
            D_t^+ \norm{Z(t)}_{\star} & = \limsup_{h \rightarrow 0^+} \frac{\norm{Z(t+h)}_{\star} - \norm{Z(t)}_{\star}}{h}                                 \\
                                     & = \lim_{h \rightarrow 0^+} \frac{\norm{Z(t) + h \cL Z(t) + hR(t)}_{\star} - \norm{Z(t)}_{\star}}{h}  \\   
                                     & \leq \lim_{h \rightarrow 0^+} \frac{\norm{Z(t) + h \cL Z(t)}_{\star} - \norm{Z(t)}_{\star}}{h} + \norm{R(t)}_{\star} \\
                                     & \leq \lim_{h \rightarrow 0^+} \frac{\norm{I + h \cL}_{2} - 1}{h} \norm{Z(t)}_{\star} + \norm{R(t)}_{\star}          \\
                                     & = \ell \norm{Z(t)}_{\star} + \norm{R(t)}_{\star},
        \end{aligned}$$
        since $\norm{\,\cdot\,}_{\star}$ is sub-multiplicative.  
    The solution of this Dini differential equation verifies (see, e.g.,~\cite[Theorem 10.1]{hairersolving})
    $$\norm{Z(t)}_{\star} \leq e^{t \ell} \norm{Z(0)}_{\star} + \int_0^t e^{(t-\tau) \ell} \norm{R(\tau)}_{\star} d\tau.$$
    From Lemma~\ref{lemma: phi functions as IVPs} with $t=h$, the claimed inequality for $k=0$ is obtained with $Z(0) = X$ and $R(t) = 0$. For $k\geq 1$, we take $Z(0) = 0$ and $R(t) = \frac{t^{k-1}}{(k-1)!} X$ which gives $Z(t) = t^k \varphi_k (t \cL) X$ by the same Lemma. The inequality above reduces to
    \[
     \norm{Z(t)}_{\star} \leq \int_0^t e^{(t-\tau) \ell} \frac{\tau^{k-1}}{(k-1)!}  d\tau  \norm{X}_{\star} = t^k \varphi_k (t \ell)  \norm{X}_{\star}.
    \]
    Combining gives the desired result.
\end{proof}

\subsection{Exponential integrators} \label{sec: exponential integrators}

Exponential integrators are a particular class of integrators well suited for parabolic partial differential equations (PDEs).
In this paper, the analysis is restricted to matrix differential equations even though the analysis and the application of exponential integrators apply to a larger class of differential equations in Banach spaces. While we start by rephrasing some results from \cite{hochbruck2005explicit} and \cite{hochbruck2010exponential}, we later derive new results specialized to our framework and necessary for the analysis of the new methods.

\subsubsection*{Exponential Euler}

Let $h>0$ be a time step and define $X_{0}^{\mathrm{E}} = X_0$. Interpolating the integral in \eqref{eq:closed form solution integral} at the initial value leads to the so-called \textit{exponential Euler} method, which iterates
\begin{align} \label{eq:exponential_euler}
    X_{n+1}^{\mathrm{E}} = e^{h \cL} X_n^{\mathrm{E}} + h \vphi{h \cL}\cG(X_n^{\mathrm{E}}).
\end{align}
According to the stability conditions \eqref{eq: equilibria conditions} discussed below, this method is also the only sensible choice for a one-stage method.
In the context of matrices, its analysis is straightforward and the following theorem states its convergence behavior. It can be proven like the related result for expontial Runge method below; we refer to Appendix~\ref{appendix: exponential euler} for a self-contained proof.
\begin{theorem}[Convergence of exponential Euler] \label{theorem: convergence of exponential Euler}
    Assume that $\cG$ is locally Lipschitz-continuous in a strip along the exact solution.
    Let $X_n^{\mathrm{E}}$ be the $n$-th iteration of the exponential Euler scheme \eqref{eq:exponential_euler}, and let $X(t_n)$ be the solution to \eqref{eq:full problem} at time $t_n = nh$.
    Then, for all $h \leq h_0$, the error verifies
    \begin{align*}
        \norm{X_n^{\mathrm{E}} - X(t_n)} \leq C \cdot h \cdot \max_{0 \leq t \leq t_n} \norm{\frac{d}{dt}\cG(X(t))},
    \end{align*}
    where the constant $C$ is explicitly available in the proof and depends on $\ell, L_{\cG}, h_0$ and $t_n$ but not on the step-size $h$ nor the number of iterations $n$.
\end{theorem}
The method is therefore order one in $h$ and robust to stiffness since the hidden constant does not depend on the Lipschitz constant of the stiff operator $\cL$.
The proof is given in Appendix \ref{appendix: exponential euler}; it is similar to the proof of the two-stage method stated below but easier. The hidden constant is relatively small (and even often $C\leq1$) for reasonable values of $L_{\cG} > 0$ and $\ell < 0$.

\subsubsection*{Exponential Runge--Kutta}
The methods are motivated by the closed form solution \eqref{eq:closed form solution integral} on which we apply quadrature formulas. The so-called \textit{exponential Runge--Kutta} methods with $s$ stages are given by
\begin{align} \label{eq:exponential_runge_kutta}
    \begin{aligned}
        G_{nj}  & = \cG(X_{nj}),                                               & \text{ for } j=1, \ldots, s \\
        X_{ni}  & = e^{c_i h \cL} X_n + h  \sum_{j=1}^s a_{ij}(h \cL) G_{nj}, & \text{ for } i=1, \ldots, s \\
        X_{n+1} & = e^{h \cL} X_n + h \sum_{i=1}^s b_i(h \cL) G_{ni}.
    \end{aligned}
\end{align}
The choice of operators $a_{ij}(h \cL), b_i(h \cL)$, and coefficients $c_i$, plays a crucial role in the behavior of the resulting method. In the next paragraphs, we derive explicit order conditions and then focus on the two-stage methods.

\subsubsection*{Explicit stiff order conditions}

A Runge-Kutta method is called explicit when $a_{ij} = 0$ for all $j>i-1$.
In the following, we discuss explicit order conditions leading to robust-to-stiffness methods.
As introduced in \cite{hochbruck2005explicit}, the first condition comes naturally by imposing the scheme to preserve the equilibria $X^{\star}$ of the autonomous problem $\dt{X}(t)~=~\mathcal L X(t)~+~\cG(X(t))$.
Assuming the scheme starts in an equilibrium, we must have $X^{\star} = X_{n} = X_{ni}$ for all $n > 0$ and $i \leq s$. 
Writing down the equations, we get the equilibrium conditions 
\begin{equation} \label{eq: equilibria conditions}
\sum_{i=1}^s b_i(h \mathcal L) = \varphi_1 (h \mathcal L), \qquad \sum_{j=1}^{i-1} a_{ij}(h \mathcal L) = c_i \varphi_1(c_i h \mathcal L), \quad 1 \leq i \leq s.
\end{equation}
The derivation of order conditions is well-known in the literature and is done for general operators in \cite{hochbruck2005explicit}, \cite{hochbruck2010exponential} and also for higher-orders in \cite{luan2013exponential} and then simplified in \cite{luan2014stiff}. Such analysis requires that the expressions $\varphi_k(h \mathcal L)$ are bounded, which we have shown in Lemma~\ref{lemma: phi functions are bounded}.

In the papers mentioned above, the analysis is typically done in a general Banach space, and the hidden constant are not given explicitly. Here, we do a similar analysis applied to finite dimensional vector spaces, which allows us to give explicit hidden constants, and therefore a better a priori control on the error. In particular, equation \eqref{eq: remainders} below is not available as such in \cite{hochbruck2005explicit} but is needed in Lemma \ref{lemma: remainders bounds} and for the rest of our analysis.

Let us start by defining the errors
\begin{equation} \label{eq: errors}
E_{ni}~=~X_{ni}-X(t_n + c_i h), \quad E_{n+1}~=~X_{n+1}-X(t_{n+1}).
\end{equation}
In the following, we are going to perform a Taylor expansion on the function $f(t) = \cG(X(t))$, where $X(t)$ is the solution to the full problem \eqref{eq:full problem}. We therefore assume that $\cG$ is sufficiently differentiable. Then, we plug the expansion in the definition of the method \eqref{eq:exponential_runge_kutta} and compare with the closed form formula:
\begin{equation} \label{eq: closed form with theta}
X(t_n + \theta h) = e^{\theta h \cL} X(t_n) + \int_0^{\theta h} e^{(\theta h - \tau) \cL} f(t_n + \tau ) d\tau.
\end{equation}
Assuming $\cG$ sufficiently differentiable, a Taylor expansion gives
\begin{equation} \label{eq: Taylor on G}
f(t_n + \tau) = \sum_{j=1}^{q} \frac{\tau^{j-1}}{(j-1)!} f^{(j-1)}(t_n) + \int_0^{\tau} \frac{(\tau - \sigma)^{q-1}}{(q-1)!} f^{(q)} (t_n + \sigma) d \sigma,
\end{equation}
and therefore we have with~\eqref{def: phi functions} that
\begin{equation} \label{eq: closed form with Taylor}
X(t_n + c_i h) = e^{c_i h \cL} X(t_n) + \sum_{j=1}^{q_i} (c_i h)^j \varphi_j(c_i h \cL) f^{(j-1)} (t_n) + \int_0^{c_i h} e^{(c_i h - \tau) \cL} \int_0^{\tau} \frac{(\tau - \sigma)^{q_i-1}}{(q_i-1)!} f^{(q_i)}(t_n + \sigma) d \sigma d\tau.
\end{equation}
Next, the exact solution agrees with the numerical schemes up to the defects:
\begin{equation} \label{eq: exact solution into numerical scheme}
\begin{aligned}
X(t_n + c_i h) &= e^{c_i h \cL} X(t_n) + h \sum_{j=1}^{i-1} a_{ij} (h \cL) f(t_n + c_j h) + \Delta_{ni}, \\
X(t_{n+1}) &= e^{h \cL} X(t_n) + h \sum_{i=1}^s b_i (h \cL) f(t_n + c_i h) + \delta_{n+1}.
\end{aligned}
\end{equation}
Now using again the Taylor expansion \eqref{eq: Taylor on G}, we get
\begin{equation} \label{eq: exact solution into numerical scheme after expansion}
\begin{aligned}
&X(t_n + c_i h) = e^{c_i h \cL} X(t_n) + h \sum_{j=1}^{i-1} a_{ij} (h \cL) \left( \sum_{k=1}^{q_i} \frac{(c_j h)^{k-1}}{(k-1)!} f^{(k-1)}(t_n) + \int_0^{c_j h} \frac{(c_j h - \sigma)^{q_i-1}}{(q_i-1)!} f^{(q_i)} (t_n + \sigma) d \sigma \right) + \Delta_{ni}, \\
&X(t_{n+1}) = e^{h \cL} X(t_n) + h \sum_{i=1}^s b_i (h \cL) \left( \sum_{k=1}^{q} \frac{(c_i h)^{k-1}}{(k-1)!} f^{(k-1)}(t_n) + \int_0^{c_i h} \frac{(c_i h - \sigma)^{q-1}}{(q-1)!} f^{(q)} (t_n + \sigma) d \sigma \right) + \delta_{n+1}.
\end{aligned}
\end{equation}
Comparing equations \eqref{eq: closed form with Taylor} and \eqref{eq: exact solution into numerical scheme after expansion} gives the following expressions for the defects:
\begin{equation} \label{eq: defects}
\begin{aligned}
&\Delta_{ni} = \sum_{j=1}^{q_i} h^j \psi_{j,i} (h \cL) f^{(j-1)}(t_n) + \Delta_{ni}^{[q_i]}, \quad &\psi_{j,i}(h \cL) = \varphi_j (c_i h \cL) c_i^j - \sum_{k=1}^{i-1} a_{ik}(h \cL) \frac{c_k^{j-1}}{(j-1)!}, \\
&\delta_{n+1} = \sum_{j=1}^q h^j \psi_j(h \cL) f^{(j-1)} (t_n) + \delta_{n+1}^{[q]}, \quad &\psi_j(h \cL) = \varphi_j (h \cL) - \sum_{k=1}^s b_k(h \cL) \frac{c_k^{j-1}}{(j-1)!},
\end{aligned}
\end{equation}
and the remainders are explicitly given by
\begin{equation} \label{eq: remainders}
\begin{split}
&\Delta_{ni}^{[q_i]} = \int_0^{c_i h} e^{(c_i h - \tau) \cL} \int_0^{\tau} \frac{(\tau - \sigma)^{q_i-1}}{(q_i-1)!} f^{(q_i)}(t_n + \sigma) d \sigma d\tau - h \sum_{j=1}^{i-1} a_{ij} (h \cL) \int_0^{c_j h} \frac{(c_j h - \sigma)^{q_i-1}}{(q_i-1)!} f^{(q_i)} (t_n + \sigma) d \sigma, \\
&\delta_{n+1}^{[q]} = \int_0^{h} e^{(h - \tau) \cL} \int_0^{\tau} \frac{(\tau - \sigma)^{q-1}}{(q-1)!} f^{(q)}(t_n + \sigma) d \sigma d\tau - h \sum_{i=1}^s b_i(h \cL) \int_0^{c_i h} \frac{(c_i h - \sigma)^{q-1}}{(q-1)!} f^{(q)} (t_n + \sigma) d \sigma.
\end{split}
\end{equation}
Finally, the general recursion for the errors \eqref{eq: errors} is given by
\begin{equation} \label{eq: error recursion}
\begin{aligned}
E_{ni} = e^{c_i h \cL} E_n + h \sum_{j=1}^{i-1} a_{ij}(h \cL) \left[ \cG(X_{nj}) - f(t_n + c_j h) \right] - \Delta_{ni}, \\
E_{n+1} = e^{h \cL} E_n + h \sum_{i=1}^s b_i(h \cL) \left[\cG(X_{ni}) - f(t_n + c_i h) \right] - \delta_{n+1}.
\end{aligned}
\end{equation}
Except for equation \eqref{eq: remainders}, the derivations above are not new and are also derived in \cite{hochbruck2005explicit}, where the analysis is then performed directly on the recursion of the error.
The new methods, however, require a different approach because of the truncation operators as we will see in Section \ref{sec:projected_exponential_methods}. We therefore need the following lemma from which we will derive our own analysis.
\begin{lemma} \label{lemma: remainders bounds}
Assume that $f(t) = \cG(X(t))$ is sufficiently differentiable. For all $h \leq h_0$, the remainders given in \eqref{eq: remainders} are bounded by
\begin{align*}
&\norm{\Delta_{ni}^{[q_i]}} \leq C_{i}^{[q_i]} \cdot h^{q_i + 1} \cdot \max_{0 \leq s \leq h} \norm{f^{(q_i)} (t_n + s)}, \\
&\norm{\delta_{n+1}^{[q]}} \leq C^{[q]} \cdot h^{q + 1} \cdot \max_{0 \leq s \leq h} \norm{f^{(q)} (t_n + s)}.
\end{align*}
The constants depend on $\ell$, the coefficients $c_i$ and $h_0$, but not on the number of steps $n$.
They are explicitly available in the proof and we have, for example for the two-stage method \eqref{eq: exponential Runge}, that
\begin{equation*}
C_1^{[1]} = 0, \quad C_2^{[1]} = (c_2 h)^2 \varphi_2(h \ell), \quad  C^{[2]} = h^3 \varphi_3(h \ell) + \frac{c_2^2 h^3}{2} \varphi_2(h \ell).
\end{equation*}
\end{lemma}

\begin{proof}
The proof is by direct calculations. For the remainder $\Delta_{ni}^{[q_i]}$:
\begin{align*}
\norm{\int_0^{c_i h} e^{(c_i h - \tau) \cL} \int_0^{\tau} \frac{(\tau - \sigma)^{q_i-1}}{(q_i-1)!} f^{(q_i)}(t_n + \sigma) d \sigma d\tau} &\leq \int_0^{c_i h} e^{(c_i h - \tau) \ell} \int_0^{\tau} \frac{(\tau - \sigma)^{q_i-1}}{(q_i-1)!} d \sigma d\tau \max_{0 \leq s \leq h} \norm{f^{(q_i)}(t_n + s)} \\
&= (c_i h)^{q_i + 1} \varphi_{q_i + 1} (h \ell) \cdot \max_{0 \leq s \leq h} \norm{f^{(q_i)}(t_n + s)}, \\
\norm{h \sum_{j=1}^{i-1} a_{ij} (h \cL) \int_0^{c_j h} \frac{(c_j h - \sigma)^{q_i-1}}{(q_i-1)!} f^{(q_i)} (t_n + \sigma) d \sigma} &\leq h \norm{\sum_{j=1}^{i-1} a_{ij}(h \cL)} \int_0^h \frac{(h - \sigma)^{q_i-1}}{(q_i-1)!} d \sigma \max_{0 \leq s \leq h} \norm{f^{(q_i)} (t_n + s)} \\
&\leq h^{q_i+1} \cdot \frac{c_i \varphi_1(c_i h \ell)}{q_i!} \cdot \max_{0 \leq s \leq h} \norm{f^{(q_i)} (t_n + s)}.
\end{align*}
In the second bound, we used the stability condition $\sum_{j=1}^{i-1} a_{ij}(h \cL) = c_i \varphi_1(c_i h \cL)$. 
Therefore, the remainder is bounded by
\begin{align*}
\norm{\Delta_{ni}^{[q_i]}}
&\leq C h^{q_i+1} \max_{0 \leq s \leq h} \norm{f^{(q_i)} (t_n + s)},
\end{align*}
for some constant $C$ independent of $h$, when $h \leq h_0$.
The remainder $\delta_{n+1}^{[q]}$ is bounded similarly:
\begin{align*}
\norm{\int_0^{h} e^{(h - \tau) \cL} \int_0^{\tau} \frac{(\tau - \sigma)^{q-1}}{(q-1)!} f^{(q)}(t_n + \sigma) d \sigma d\tau} &\leq \int_0^{h} e^{(h - \tau) \ell} \int_0^{\tau} \frac{(\tau - \sigma)^{q-1}}{(q-1)!}  d \sigma d\tau \cdot \max_{0 \leq s \leq h} \norm{f^{(q)}(t_n + s)}, \\
&= h^{q+1} \varphi_{q+1}(h \ell) \cdot \max_{0 \leq s \leq h} \norm{f^{(q)}(t_n + s)} \\
\norm{h \sum_{i=1}^s b_i(h \cL) \int_0^{c_i h} \frac{(c_i h - \sigma)^{q-1}}{(q-1)!} f^{(q)} (t_n + \sigma) d \sigma} &\leq h^{q+1} \cdot \frac{\varphi_1(h \ell)}{q!} \cdot \max_{0 \leq s \leq h} \norm{f^{(q)} (t_n + s)},
\end{align*}
where we used the stability condition $\sum_{i=1}^s b_i(h \mathcal L) = \varphi_1 (h \mathcal L)$. The conclusion follows immediately.
\end{proof}

The order conditions are naturally obtained in expression \eqref{eq: defects} by imposing $\psi_{ij} = 0$ and $\psi_k = 0$ for all indexes $i,j,k$.

\subsubsection*{Exponential Runge}

The second order conditions are obtained by imposing 
\begin{equation} \label{eq: order condition two stages}
\psi_{1,2}(h \cL) \overset{!}{=} 0, \quad \text{and} \quad \psi_1 (h \cL) \overset{!}{=} 0 \overset{!}{=} \psi_2 (h \cL).
\end{equation}
Combining with the equilibrium conditions \eqref{eq: equilibria conditions}, we obtain the following family of explicit methods for $c_2 \neq 0$, called the exponential Runge methods, in a Butcher table:
\begin{equation} \label{eq: exponential Runge table}
\begin{array}{c|cc}
0 \\
c_2 & c_2 \varphi_{1,2}  \\
\hline
& \varphi_1 - \frac{1}{c_2} \varphi_2  & \frac{1}{c_2} \varphi_2
\end{array}
\end{equation}
Written more explicitly, this gives us the two-stage method
\begin{equation} \label{eq: exponential Runge}
\left\{
\begin{aligned}
    & X_{n2}^{\mathrm{R}} = e^{c_2 h \cL} X_n^{\mathrm{R}} + c_2 h \, \varphi_1 \left(c_2 h \cL \right) \cG(X_n^{\mathrm{R}}), \\
    & X_{n+1}^{\mathrm{R}}     = e^{h \cL} X_n^{\mathrm{R}} + h \varphi_1(h \cL) \cG(X_n^{\mathrm{R}}) + h \varphi_2(h \cL) \left( \cG(X_{n2}^{\mathrm{R}}) - \cG(X_n^{\mathrm{R}}) \right).
\end{aligned}
\right.
\end{equation}
The following theorem proves the convergence of these methods, and its proof gives information on the hidden constants.

\begin{theorem}[Convergence of exponential Runge] \label{theorem: convergence of exponential Runge}
    Assume that $\cG$ is locally Lipschitz-continuous in a strip along the exact solution, and sufficiently differentiable.
    Let $X_n^{\mathrm{R}}$ be the $n$-th iteration of the exponential Runge scheme \eqref{eq: exponential Runge} started at $X_0 = X(0)$, and let $X(t_n)$ be the solution to \eqref{eq:full problem} at time $t_n = nh$.
    Then, for all $h \leq h_0$, the error verifies
    \begin{align*}
        \norm{X_n^{\mathrm{R}} - X(t_n)} \leq C \cdot h^2 \cdot \left( \max_{0 \leq t \leq t_n} \norm{\frac{d}{dt}\cG(X(t))} + \max_{0 \leq t \leq t_n} \norm{\frac{d^2}{dt^2}\cG(X(t))} \right),
    \end{align*}
    where the constant $C$ depends on $\ell, L_{\cG}, h_0$ and $t_n$ but not on $h$ nor $n$.
\end{theorem}

\begin{proof}
From \eqref{eq: error recursion} and taking the norm, the error verifies
\begin{equation} \label{eq: runge error recursion}
\begin{aligned}
\norm{E_{n+1}} &= \norm{e^{h \cL} E_n + h \sum_{i=1}^2 b_i(h \cL) \left[ \cG(X_{ni}) - f(t_n + c_i h) \right] + \delta_{n+1}} \\
&\leq e^{h \ell} \norm{E_n} + h \sum_{i=1}^2 \norm{b_i (h \cL)} \cdot L_{\cG} \norm{E_{ni}} + \norm{\delta_{n+1}}.
\end{aligned}
\end{equation}
By definition of the exponential Runge method, we have that
$$\norm{E_{n1}} = \norm{E_n} + \norm{\Delta_{n1}}, \quad \norm{E_{n2}} \leq e^{c_2 h \ell} \norm{E_n} + c_2 h \varphi_1(c_2 h \ell) \norm{E_{n1}} + \norm{\Delta_{n2}},$$
and from the order conditions \eqref{eq: order condition two stages} and Lemma \ref{lemma: remainders bounds}, it follows that
\begin{equation} \label{eq: remainders for exponential Runge}
\begin{aligned}
\norm{\Delta_{n1}} = \norm{\Delta_{n1}^{[1]}} = 0, \qquad 
\norm{\Delta_{n2}} = \norm{\Delta_{n2}^{[1]}} 
\leq (c_2 h)^2 \varphi_2(h \ell) \cdot \max_{0 \leq s \leq h} \norm{f'(t_n + s)}, \\
\text{and } \quad \norm{\delta_{n+1}} = \norm{\delta_{n+1}^{[2]}}
\leq \left[h^3 \varphi_3(h \ell) + \frac{c_2^2 h^3}{2} \varphi_2(h \ell) \right] \cdot \max_{0 \leq s \leq h} \norm{f''(t_n + s)}.
\end{aligned}
\end{equation}
Next, we define the following quantities:
\begin{equation} \label{eq: beta1 and beta2}
\begin{aligned}
\beta_1 = (c_2 h)^2 \varphi_2(h \ell) \cdot \max_{0 \leq s \leq t_{n+1}} \norm{f'(s)}, \qquad
\beta_2 = \left[h^2 \varphi_3(h \ell) + \frac{c_2^2 h^2}{2} \varphi_2(h \ell) \right] \cdot \max_{0 \leq s \leq t_{n+1}} \norm{f''(s)}.
\end{aligned}
\end{equation}
When $h \leq h_0 \leq c_2$, direct computations give that
\begin{equation} \label{eq: bounds for bi}
\norm{b_1(h \cL)} \leq \varphi_1(h \ell) - \tfrac{1}{c_2} \varphi_2(h \ell), \qquad \norm{b_2(h \cL)} \leq \tfrac{1}{c_2} \varphi_2(h \ell).
\end{equation}
Now inserting \eqref{eq: remainders for exponential Runge} into \eqref{eq: runge error recursion} and simplifying with \eqref{eq: beta1 and beta2} and \eqref{eq: bounds for bi}, we obtain the recursion
\begin{equation*}
\norm{E_{n+1}} \leq (1 + h L_*) \norm{E_n} + h (\beta_1 + \beta_2),
\end{equation*}
where $L_* = \ell + \varphi_1(h \ell) L_{\cG} + O(h)$, and the hidden constant in $O(h)$ depends only on $\ell$ and $L_{\cG}$ (multiplied by $\varphi$-functions).
Solving the recursion leads to
\begin{equation} \label{eq: solving recursion}
\norm{E_n} \leq (1+h L_*)^n \norm{E_0} + h \sum_{k=0}^{n-1} (1+h L_*)^{k} (\beta_1 + \beta_2),
\end{equation}
Since we start from the initial value $X_0 = X(0)$, it holds  $\norm{E_0} = 0$. With the inequalitiy
$$
h\sum_{k=0}^{n-1} (1 + h L_*)^k \leq \frac{e^{nh L_*} - 1}{L_*},$$
we conclude that
\begin{equation*}
\norm{X_n^{\mathrm{R}} - X(t_n)} \leq \frac{e^{n h L_*} - 1}{L_*} \cdot (\beta_1 + \beta_2).
\end{equation*}
For  all $h \leq h_0$, the factor above reduces to a constant $C(\ell, L_{\cG}, h_0)$ as stated in the theorem.
\end{proof}

\begin{remark} \label{remark: non strict Runge}
The exponential Runge methods~\eqref{eq: exponential Runge table} require the functions $\varphi_1$ and $\varphi_2$. It is tempting to consider the following family of methods that only requires $\varphi_1$:
\begin{equation*} \label{eq: exponential Runge not order 2}
\begin{array}{c|cc}
0 \\
c_2 & c_2 \varphi_{1,2}  \\
\hline
& (1 - \frac{1}{2 c_2}) \varphi_1 & \frac{1}{2 c_2} \varphi_1
\end{array}
\end{equation*}
These methods have a classical order of two but their stiff order is strictly smaller than 2; see \cite{hochbruck2005explicit} for more details. Nevertheless, since their evaluation per step is much cheaper than the exponential Runge methods, they might be useful in some practical applications. In the following, we will refer to these methods as non-strict exponential Runge methods.
\end{remark}

\subsection{Dynamical low-rank approximation}

Let $\Mr$ denote the set of $m \times n$ matrices of rank $r$. It is a smoothly embedded submanifold in $\Rmn$ and its Riemannian geometry can be used effectively in numerical algorithms; see~\cite{uschmajew2020geometric} for an overview.

Suppose the solution of the full problem \eqref{eq:full problem} admits a good low-rank approximation. In that case, one can apply the dynamical low-rank approximation (DLRA) and solve a reduced problem instead. The technique is extensively studied in \cite{koch2007dynamical}. The corresponding DLRA of rank $r$ is
\begin{align*}
    \dt{Y}(t) = \proj{Y(t)}{\cL Y(t) + \cG (Y(t))}, \quad Y(0) = Y_0 \in \Mr,
\end{align*}
where $\mathcal{P}_{Y}$ is the $l_2$-orthogonal projection onto the tangent space $\mathcal T_{Y} \Mr$ of $\Mr$ at $Y \in \Mr$, 
$$
 \mathcal{P}_{Y}(Z) = UU^T Z + Z VV^T - UU^T Z VV^T \quad \text{with $Y=U\Sigma V^T$ a compact SVD},
$$
and $Y_0 \approx X_0$ is a (near) best rank approximation of $X_0$. Since $\cL Y \in \mathcal T_Y \Mr$, the equation simplifies to
\begin{equation} \label{eq:DLRA problem after linearity}
    \dt{Y}(t) = \cL Y(t) + \proj{Y(t)}{\cG(Y(t))}, \quad Y(0) = Y_0.
\end{equation}
\begin{remark}
As long as the solution of \eqref{eq:DLRA problem after linearity} exists, it belongs to the manifold $\Mr$.
\end{remark}
Recall that $\phi_{\mathcal F}^t(X_0)$ denotes the flow of the full problem \eqref{eq:full problem} at time $t$ with initial value $X_0$.
Similarly, we denote by $\phi_{P \cF}^t(Y_0)$ the flow of the projected problem \eqref{eq:DLRA problem after linearity} with initial value $Y_0 \in \Mr$.
The following assumption is standard when studying the approximation error made by DLRA.
\begin{assumption}[DLRA assumption] \label{ass:DLRA assumptions}
    The function $\cG$ maps almost to the tangent bundle of $\Mr$:
    $$\norm{\cG(Y)-\mathcal{P}_Y \cG(Y)} \leq \varepsilon_r,$$
    for all $Y$ in a neighborhood of the solution to the DLRA.
\end{assumption}
In terms of the total vector field, $\mathcal F = \cL + \cG$, the assumption above is equivalent to $\norm{\mathcal F(Y)-\mathcal{P}_Y \mathcal F(Y)}~\leq~\varepsilon_r$ since $\cL$ maps to the tangent space. Using this property, we get immediately the following approximation result of DLRA (see~\cite{kieri2016discretized} for a proof):

\begin{theorem}[Accuracy of DLRA] \label{thm:DLRA theoretical error}
    Assuming that $\mathcal F = \mathcal L + \mathcal G$ is one-sided Lipschitz with constant $\ellF$, and Assumption \ref{ass:DLRA assumptions}, the error made by the dynamical low-rank approximation of rank $r$ verifies
    \begin{align*}
        \norm{\phi_{\cF}^t(X_0) - \phi_{P \cF}^t(Y_0)} \leq e^{t \ellF} \norm{X_0 - Y_0} + t \varphi_1(t \ell_\cF) \, \varepsilon_r ,
    \end{align*}
    where $\mathcal F$ represents the vector field of the full problem \eqref{eq:full problem}, and $P \mathcal F$ that of the projected problem~\eqref{eq:DLRA problem after linearity}.
\end{theorem}
Recall the convenient bound $\ellF \leq \ell_{\cL}+ L_{\cG}$ from Remark~\ref{remark: bound ellF sum}, showing that $\ellF$ can be negative for small non-linear perturbations $\mathcal G$ with diffusive $\cL$.

\subsection{Contributions and outline}

\begin{figure}[ht]
    \centering
    \includegraphics[width=0.7\textwidth]{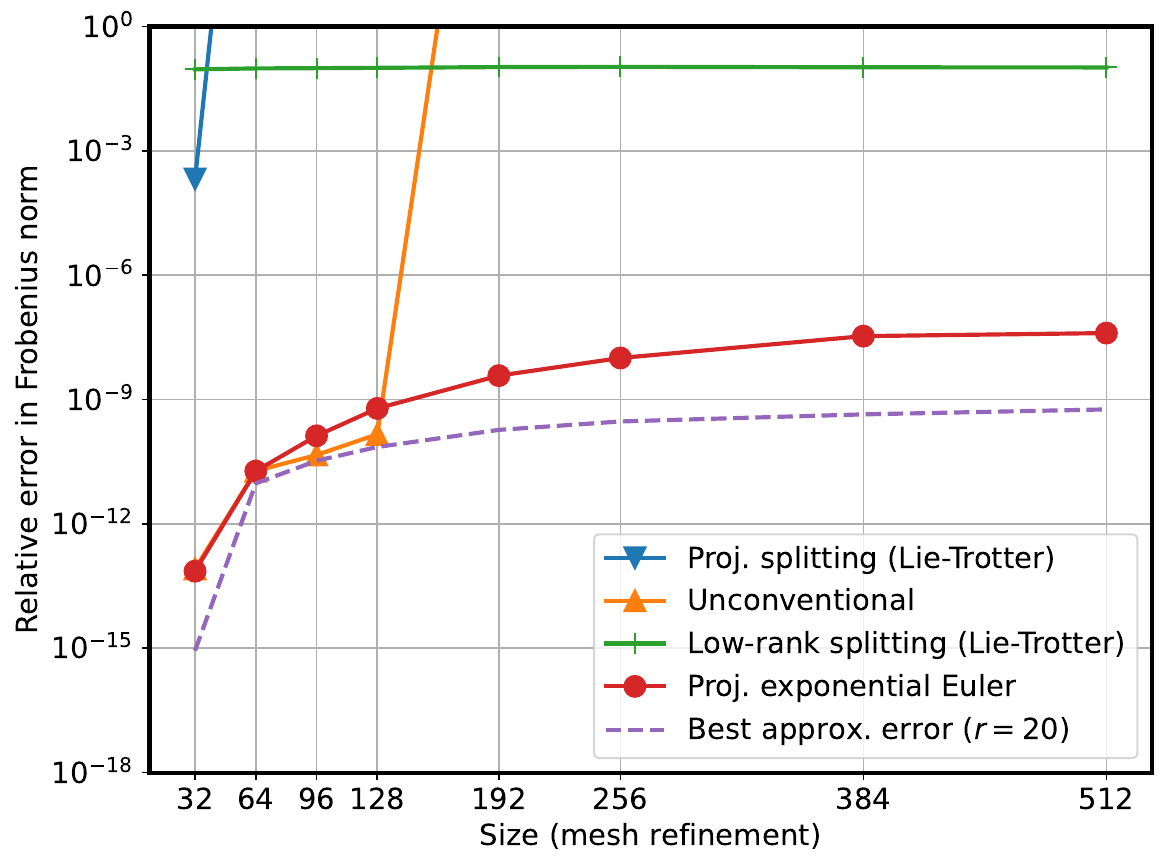}
    \caption{Error of several methods under mesh refinements for solving the heat equation (stiff) that can be formulated as differential Lyapunov equation with constant source $C$; see Section \ref{sec:diff Lyap}. Step size is $h = 0.01$.}
    \label{fig:motivation}
\end{figure}

Even though they do not impose a particular structure to the problem, the projector-splitting algorithm \cite{lubich2014projector} and the unconventional integrator \cite{ceruti2022unconventional} are not proven to be robust to stiff problems. 
In fact, Figure~\ref{fig:motivation} shows that the two methods are numerically not robust to stiffness in general, even when a very small step size ($10^{-5}$) is used for the sub-steps.
In addition, the projector-splitting algorithm contains a backward sub-step, making it not suitable for parabolic problems.
More specialized methods are needed for stiff problems.
Assuming the problem has a Sylvester-like structure, the low-rank splitting proposed in \cite{ostermann2019convergence} is shown theoretically robust to stiffness but the method is not very useful in practice because of its large hidden constant.
When the step size is $h \approx 0.01$, we can see in Figure~\ref{fig:motivation} that its error is unreasonably large, though robust under mesh refinements.
This paper's contribution is to fill the gap by proposing new methods of integration of DLRA for Sylvester-like structures that are robust to stiffness and have small hidden error constants. 
As a motivation, Figure~\ref{fig:motivation} shows that one of these new methods is accurate and robust under mesh refinements. 
More numerical experiments are detailed in Section \ref{sec:numerical_experiments}.

The paper is organized as follows. Section \ref{sec:projected_exponential_methods} defines the new methods and derives theoretical guarantees. Section \ref{sec: efficient implementation} describes an efficient implementation of the new methods, making this new method a practical integrator. Section \ref{sec:numerical_experiments} shows numerical experiments of the method and verifies the theoretical guarantees derived in the paper. Section \ref{sec:conclusion} concludes and discusses potential future work directions.

\section{Projected exponential methods} \label{sec:projected_exponential_methods}

The definition of the (explicit) \textit{projected exponential Runge-Kutta} (PERK) methods is naturally obtained from~\ref{eq:exponential_runge_kutta}. 
We state them here for a method with $s$ stages as
\begin{equation} \label{eq: projected exponential methods}
    \begin{aligned}
        G_{nj}  & = \proj{Y_{nj}}{\cG (Y_{nj})}, && \text{ for } j=1, \ldots, s                                    \\
        Y_{ni}  & = \Ret{e^{c_i h \cL} Y_n + h \sum_{j=1}^{i-1} a_{ij} (h \cL) G_{nj}}, && \text{ for } i = 1, \ldots, s \\
        Y_{n+1} & = \Ret{e^{h \cL} Y_n + h \sum_{i=1}^s b_i (h \cL) G_{ni}},
    \end{aligned}
\end{equation}
where $\Tr \colon \Rmn \rightarrow \Mr$ represents taking a best rank $r$ approximation of a matrix in Frobenius norm. Recall that this can be computed by truncating the SVD of the given matrix to its $r$ dominant terms. We therefore call this operator the truncation operator.
It is natural to consider using the same operators $a_{ij}(h \cL), b_i(h \cL)$ and coefficients $c_i$ as derived in the classical theory. Indeed, the projected methods coincide with the classical methods when $r=\min\{m, n\}$. Also, the same order conditions will appear in the expansion of the defects, as we will see later in the analysis.

The methods~\eqref{eq: projected exponential methods} are computationally attractive only when the internal stages remain of sufficiently low rank and when this is exploited in implementation. Thanks to the tangent space projection, $\rank(G_{nj})~\leq~2r$ and the growth in rank comes mainly from the terms $\sum_{j=1}^s a_{ij} (h \cL) G_{nj}$ and $\sum_{i=1}^s b_i (h \cL) G_{ni}$. The computational and memory requirements will therefore grow with the number of internal stages $s$. However, compared to classical Runge--Kutta methods, the scaling is higher than linear in $s$ due to the numerical linear algebra of low-rank calculations. We therefore limit our analysis to the one-stage and two-stage methods stated below.

\subsection{Preliminary results} \label{sec: preliminary results}

We start with Lemma~\ref{lemma:retraction inequality}  on the stability of the rank $r$ truncation operator $\Tr$ under perturbation. While it is well known that the Lipschitz constant of $\Tr$ at $X \in \Mr$ behaves like $\sigma_r^{-1}(X)$, which would lead to undesirable bounds when $\sigma_r(X)$ is small, the operator $\Tr$ is stable if we study the perturbation in an additive way. This result is similar to other results in the literature (see, e.g., \cite[Lemma 4.1]{hackbusch2016new}) but we include the proof for completeness.

\begin{lemma}[Truncation inequality] \label{lemma:retraction inequality}
    Denote by $\Tr$ the truncation to rank $r$. Then, for any two matrices $A$ and $E$, we have
    $$\norm{\Ret{A+E} - A} \leq \norm{\Ret{A}-A} + 2 \norm{E}.$$
\end{lemma}

\begin{proof}
    The triangular inequality gives
    $$\norm{\Ret{A + E} - A} = \norm{\Ret{A + E} - (A+E) + (A+E) - A} \leq \norm{\Ret{A + E} - (A+E)} + \norm{E},$$
    and the conclusion follows from the definition of the best rank $r$ approximation,
    \begin{align*}
    \norm{\Ret{A + E} - (A+E)} &= \min_{\rank(X) \leq r} \norm{X - (A+E)} 
    = \min_{\rank(X) \leq r} \norm{X - A + A - (A+E)} \\
    &\leq \min_{\rank(X) \leq r} \norm{X - A}_F +\norm{E} 
    = \norm{\Ret{A} - A} + \norm{E}. \qedhere
    \end{align*}
\end{proof}

The factor two appearing in the truncation inequality above is inevitable and brings an extra challenge to the convergence proofs.
Indeed, we cannot immediately express the numerical flow $Y_n$ as a perturbation of the solution $X(t_n)$ nor $Y(t_n)$, since we would have a factor $2^n$ appearing when solving the recursion. The issue is bypassed by comparing instead the numerical flow with the DLRA started at the previous numerical solution $\phi_{P \cF}^h(Y_n)$. This way, the exponential term is treated exactly and the undesirable factor $2^n$ does not appear. Similarly to the derivations made in Section \ref{sec: exponential integrators}, we derive a stable\footnote{Stable w.r.t.~small singular values since the expansion does not involve derivatives of the projection.} expansion of the projected field flow.
The closed form formula \eqref{eq:closed form solution integral} applied to the projected field \eqref{eq:DLRA problem after linearity} gives for $i=1, \ldots, s$
\begin{equation} \label{eq: closed form projected solution}
\phi_{P \cF}^{c_i h} (Y_n) = e^{c_i h \cL} Y_n + \int_0^{c_i h} e^{(c_i h - \tau) \cL} \proj{\phi_{P \cF}^{\tau}(Y_n)}{\cG(\phi_{P \cF}^{\tau}(Y_n))} d \tau.
\end{equation}
The next step is to use $f(t) = \cG(X(t))$ and decompose the projection as
$$\proj{\phi_{P \cF}^{\tau}(Y_n)}{\cG(\phi_{P \cF}^{\tau}(Y_n))} = f(t_n + \tau) + \proj{\phi_{P \cF}^{\tau}(Y_n)}{\cG(\phi_{P \cF}^{\tau}(Y_n))} - f(t_n+\tau).$$
Therefore, we obtain that
$$\phi_{P \cF}^{c_i h} (Y_n) = e^{c_i h \cL} Y_n + \int_0^{c_i h} e^{(c_i h - \tau) \cL}  f(t_n + \tau) d \tau + \Omega_{ni} + \Lambda_{ni},$$
where we have defined the two modelling defects
\begin{align} 
\Omega_{ni} &= \int_0^{c_i h} e^{(c_i h - \tau) \cL} \left[ \proj{\phi_{P \cF}^{\tau}(Y_n)}{\cG(\phi_{P \cF}^{\tau}(Y_n))} - \cG(\phi_{P \cF}^{\tau}(Y_n)) \right] d\tau, \label{eq: modelling defect 1} \\
\Lambda_{ni} &= \int_0^{c_i h} e^{(c_i h - \tau) \cL} \left[ \cG(\phi_{P \cF}^{\tau}(Y_n)) - f(t_n + \tau) \right] d\tau.
\label{eq: modelling defect 2}
\end{align}
By definition of $\Delta_{ni}$ in \eqref{eq: exact solution into numerical scheme} we already have
$$\int_0^{c_i h} e^{(c_i h - \tau) \cL}  f(t_n + \tau) d \tau = \sum_{j=1}^{i-1} a_{ij}(h \cL) f(t_n + c_j h) + \Delta_{ni},$$
and we therefore obtain the following expressions for the projected field flow:
\begin{equation} \label{eq: projected field expansions}
\begin{split}
\phi_{P\cF}^{c_i h}(Y_n) = e^{c_i h \cL} Y_n + h \sum_{j=1}^{i-1} a_{ij}(h \cL) f(t_n + c_j h) + \Delta_{ni} + \Omega_{ni} + \Lambda_{ni}, \\
\phi_{P\cF}^h(Y_n) = e^{h \cL} Y_n + h \sum_{i=1}^s b_i(h \cL) f(t_n + c_i h) + \delta_{n+1} + \omega_{n+1} + \lambda_{n+1},
\end{split}
\end{equation}
where the second formula was derived in a similar manner.
The next lemma bounds the modelling defects, and will be essential in the proof of convergence of the new methods.
\begin{lemma}[Stable expansion of the projected flow] \label{lemma: expansion of the projected flow}
Assume the function $\cG$ Lipschitz-continuous, and the function $f(t) = \cG(X(t))$ sufficiently differentiable. Then, under Assumption~\ref{ass:DLRA assumptions} the projected field flow verifies the expansions \eqref{eq: projected field expansions} with terms $\Delta_{ni}$ and $\delta_{n+1}$ defined in \eqref{eq: defects}, and modelling defects bounded by
$$\norm{\Omega_{ni}} \leq c_i h \varphi_1(c_i h \ell) \varepsilon_r, \quad \norm{\omega_{n+1}} \leq h \varphi_1(h \ell) \varepsilon_r,$$
$$\norm{\Lambda_{ni}} \leq L_{\cG} \cdot c_i h \varphi_1(c_i h \ell) \cdot (1 + O(c_i h)) \cdot \norm{Y_n - X(t_n)} + O((c_i h)^2) \cdot \varepsilon_r,$$
$$\norm{\lambda_{n+1}} \leq L_{\cG} \cdot h \varphi_1(h \ell) \cdot (1 + O(h)) \cdot \norm{Y_n - X(t_n)} + O(h^2) \cdot \varepsilon_r.$$
\end{lemma}

\begin{proof}
Bounding the first modelling defect \eqref{eq: modelling defect 1} requires Assumption \ref{ass:DLRA assumptions}, and a direct computation gives
\begin{align*}
\norm{\Omega_{ni}} &\leq \int_0^{c_i h} e^{(c_i h - \tau) \ell} \norm{\proj{\phi_{P \cF}^{\tau}(Y_n)}{\cG(\phi_{P \cF}^{\tau}(Y_n))} - \cG(\phi_{P \cF}^{\tau}(Y_n))} d\tau \leq c_i h \varphi_1(c_i h \ell) \varepsilon_r.
\end{align*}
The statement for $\omega_{n+1}$ holds with a step size $h$ instead of $c_i h$.

For the second modelling defect \eqref{eq: modelling defect 2}, we need to use the Lipschitz continuity of $\cG$, and then Theorem~\ref{thm:DLRA theoretical error} to get that
\begin{align*}
\norm{\Lambda_{ni}} &\leq \int_0^{c_i h} e^{(c_i h - \tau) \ell} \cdot L_{\cG} \cdot \norm{\phi_{P \cF}^{\tau}(Y_n) - \phi_{\cF}^{\tau}(X(t_n))} d \tau \\
&\leq \int_0^{c_i h} e^{(c_i h - \tau) \ell} \cdot L_{\cG} \cdot \left[ e^{\tau \ellF} \norm{Y_n - X(t_n)} + \tau \varphi_1( \tau \ellF ) \varepsilon_r \right] d \tau.
\end{align*}
Computing the integrals and then simplifying with Taylor expansion gives
\begin{align*}
\int_0^{c_i h} e^{(c_i h - \tau) \ell} e^{\tau \ellF} d\tau = \frac{e^{c_i h \ell} - e^{c_i h \ellF}}{\ell - \ellF} = c_i h \varphi_1(c_i h \ell) \cdot (1 + O(c_i h)), \\
\int_0^{c_i h} e^{(c_i h - \tau) \ell} \tau \varphi_1(\tau \ellF) d \tau = \frac{\ell (1 - e^{c_i h \ellF}) + \ellF (e^{c_i h \ell} -1)}{\ell\, \ellF (\ell - \ellF)} = O((c_i h)^2).
\end{align*}
Since $\ellF \leq \ell + L_{\cG}$, the hidden constants depend only on $\ell$ and $L_{\cG}$ and they do not require a more explicit treatment in the analysis.
\end{proof}

\subsection{Projected exponential Euler}

Similarly to the exponential Euler scheme \eqref{eq:exponential_euler}, we define one step of the \textit{projected exponential Euler} scheme:
\begin{align} \label{eq:projected exponential euler}
    Y_{n+1}^{\mathrm{E}} = \Ret{ e^{h \cL} Y_n^{\mathrm{E}} + h \varphi_1 (h \cL) \proj{Y_n^{\mathrm{E}}}{\cG (Y_n^{\mathrm{E}})} }, \quad Y_n^{\mathrm{E}} \in \Mr.
\end{align}
It corresponds to \eqref{eq: projected exponential methods} where $b_1 = \varphi_1$ and $c_1 = 1$.
The computational cost of~\eqref{eq:projected exponential euler} mainly comes from evaluating the $\varphi$-functions. 
In practice, these operators are never formed explicitly but their action on $Y_n^{\mathrm{E}} $ and $\proj{Y_n^{\mathrm{E}} }{\cG (Y_n^{\mathrm{E}} )}$ is computed directly. In particular, since $Y_n^{\mathrm{E}} $ and $\proj{Y_n^{\mathrm{E}} }{\cG (Y_n^{\mathrm{E}} )}$ have rank at most $r$ and $2r$, resp., we can use Krylov techniques to approximate the action of these operators onto a tall matrix with at most $2r$ columns. A practical implementation that exploits even more structure will be detailed in~Section~\ref{sec: efficient implementation}. The following theorem states the convergence of the method.

\begin{theorem}[Projected exponential Euler convergence] \label{theorem: projected exponential Euler convergence}

Assume that $\cG$ is locally Lipschitz-continuous in a strip along the solution $X(t)$ to~\eqref{eq:full problem}.
Let $Y_n^{\mathrm{E}}$ be the $n$-th iteration of the projected exponential Euler scheme~\eqref{eq:projected exponential euler}, and let $X(t_n)$ be the solution to~\eqref{eq:full problem} at time $t_n = nh$.
Then, under Assumption~\ref{ass:DLRA assumptions} for $h \leq h_0$, the error verifies
\begin{align*}
    \norm{Y_n^{\mathrm{E}} - X(t_n)} \leq e^{t_n L_*} \cdot \norm{Y_0 - X_0} + C_1 \cdot h \cdot \max_{0 \leq t \leq t_n} \norm{\frac{d}{dt} \cG (X(t))} + C_2 \cdot \varepsilon_r,
\end{align*}
where $L_* = \ellF + 4 \varphi_1(h \ell) L_{\cG} + O(h)$.
The constants $C_1$ and $C_2$ are explicitly available in the proof and depend only on $\ell, L_{\cG}, h_0$ and $t_n$ but not on $h, n$ nor the curvature of the manifold $\Mr$.

\end{theorem}

\begin{proof}
We omit the notation $\cdot^{\mathrm{E}}$ in the proof.
Let us define the error $E_{n+1} = Y_{n+1} - X(t_{n+1})$; we are going to find a recursion on the error as in the proof of Theorem~\ref{theorem: convergence of exponential Runge} but with extra terms due to the modelling error induced by the DLRA.
The factor two appearing in the retraction inequality from Lemma~\ref{lemma:retraction inequality} forces us to treat exactly the exponential term $e^{h \cL} Y_n$, so we need to work with the DLRA flow started at the numerical solution:
\begin{equation} \label{eq: proj exp euler retraction inequality}
\norm{Y_{n+1} - X(t_{n+1})} \leq \norm{\Tr(\phi_{P \cF}^h (Y_n) + K_{n+1}) - \phi_{P \cF}^h (Y_n) - M_{n+1}} \leq 2 \norm{K_{n+1}} + \norm{M_{n+1}},
\end{equation}
where we have defined 
$$M_{n+1} = \phi_{\cF}^h (X(t_n)) - \phi_{P \cF}^h (Y_n),$$
$$K_{n+1} = e^{h \cL} Y_n + h \varphi_1(h \cL) \proj{Y_n}{\cG(Y_n)} - \phi_{P \cF}^h (Y_n).$$
The first difference term is easily bounded with Theorem~\ref{thm:DLRA theoretical error}:
\begin{equation} \label{eq: proj exp euler bound M}
\norm{M_{n+1}} \leq e^{h \ellF} \norm{X(t_n) - Y_n} + h \varphi_1(h \ellF) \varepsilon_r.
\end{equation}
The subtlety is in the second difference term $K_{n+1}$, for which we require the stable expansion of the projected flow derived in Section~\ref{sec: preliminary results}. Substituting $\phi_{P \cF}^h (Y_n)$ with~\eqref{eq: projected field expansions} and applying Lemma~\ref{lemma: expansion of the projected flow} with $b_1 = \varphi_1$ gives that
\begin{equation} \label{eq: proj exp euler bound K}
\begin{aligned}
\norm{K_{n+1}} &\leq h \varphi_1(h \ell) \cdot \| \proj{Y_n}{\cG(Y_n)} - f(t_n) \| + \norm{\delta_{n+1}} + \norm{\omega_{n+1}} + \norm{\lambda_{n+1}} \\
&\leq h \varphi_1(h \ell) \cdot (2 + O(h)) \cdot L_{\cG} \norm{E_n} + h^2 \varphi_2(h \ell) \cdot \max_{0 \leq s \leq h} \norm{f'(t_n+s)} + 2 h \varphi_1(h \ell) \varepsilon_r + O(h^2) \varepsilon_r.
\end{aligned}
\end{equation}
Inserting~\eqref{eq: proj exp euler bound K} and~\eqref{eq: proj exp euler bound M} into~\eqref{eq: proj exp euler retraction inequality} gives the following recursion of the error $E_n^Y = Y_n - X(t_n)$:
\begin{align*}
\norm{E^Y_{n+1}} &\leq \left( e^{h \ellF} + h \varphi_1(h \ell) L_{\cG} \cdot (4 + O(h)) \right) \norm{E_n^Y} \\
&\quad + 2 h^2 \varphi_2(h \ell) \cdot \max_{0 \leq s \leq t_n} \norm{f'(s)} + ( 4 h \varphi_1(h \ell) + h \varphi_1(h \ellF) ) \varepsilon_r + O(h^2) \varepsilon_r \\
&\leq \left( 1 + h L_* \right) \norm{E_n^Y} + h \left(\beta + \gamma \right),
\end{align*}
where $L_* = \ellF + 4 \varphi_1(h \ell) L_{\cG} + O(h)$ and where we have defined
\begin{align*}
\beta = 2 h \varphi_2(h \ell) \cdot \max_{0 \leq s \leq t_n} \norm{f'(s)}, \quad \gamma =  \left( 4 \varphi_1(h \ell) + \varphi_1(h \ellF) + O(h) \right) \cdot \varepsilon_r.
\end{align*}
Finally, solving the recursion as we have already done in \eqref{eq: solving recursion} leads to
$$\norm{E_n^Y} \leq e^{nh L_*} \norm{Y_0 - X_0} + \frac{e^{n h L_*} - 1}{L_*} \cdot ( \beta + \gamma ),$$
and by tracking the hidden constants, we see that they depend only on $\ell, L_{\cG}$ and $h_0$.
\end{proof}

\subsection{Projected exponential Runge}

As we can see from Lemma \ref{lemma: expansion of the projected flow}, the remainder terms $\Delta_{ni}$ and $\delta_{n+1}$ are unchanged, and therefore the Butcher table \ref{eq: exponential Runge table} for the coefficients is still valid.
The \textit{projected exponential Runge} scheme iterates
\begin{equation} \label{eq: proj exp Runge}
\begin{aligned}
    Y_{n2}^{\mathrm{R}} & = \Ret{e^{c_2 h \cL} Y_n^{\mathrm{R}} + c_2 h \varphi_1 \left( c_2 h \cL \right) \proj{Y_n^{\mathrm{R}}}{\cG(Y_n^{\mathrm{R}})}}, \\
    Y_{n+1}^{\mathrm{R}}     & = \Ret{e^{h \cL} Y_n^{\mathrm{R}} + h \varphi_1 \left( h \cL \right) \proj{Y_n^{\mathrm{R}}}{\cG(Y_n^{\mathrm{R}})} + h \varphi_2(h \cL) \left( \proj{Y_{n2}^{\mathrm{R}}}{\cG(Y_{n2}^{\mathrm{R}})} - \proj{Y_n^{\mathrm{R}}}{\cG(Y_n^{\mathrm{R}})} \right)}.
\end{aligned}
\end{equation}
This method has a good balance of numerical accuracy (order) and computational efficiency (cost per step).
Even though the scheme~\eqref{eq: proj exp Runge} requires evaluating a $\varphi_2$ function, its computation is still manageable as we will show further in the implementation section. The main convergence theorem relies heavily on Lemma~\ref{lemma: expansion of the projected flow}, but is similar in spirit to the proof of Theorem~\ref{theorem: convergence of exponential Runge}.

\begin{theorem}[Projected exponential Runge convergence]
Assume $\cG$ sufficiently many times Fréchet differentiable in a strip along the exact solution and Assumption \ref{ass:DLRA assumptions} satisfied. Let $Y_n^{\mathrm{R}}$ be the $n$-th step of the projected exponential Runge scheme started at $Y_0 \in \Mr$ with $c_2 > 0$, and $X(t_n)$ be the solution to the full problem \eqref{eq:full problem}. Then, under Assumption~\ref{ass:DLRA assumptions} for $h \leq h_0$, the error verifies
\begin{align*}
    \norm{Y_n^{\mathrm{R}} - X(t_n)} &\leq e^{t_n L_*} \norm{Y_0 - X_0} + C_1 \cdot h^2 \cdot \left( \max_{0 \leq t \leq t_n} \norm{\frac{d}{dt} \cG(X(t))} + \max_{0 \leq t \leq t_n} \norm{\frac{d^2}{dt^2} \cG(X(t))} \right) + C_2 \cdot \varepsilon_r,
\end{align*}
where $L_* = \ell_F + 4 \varphi_1(h \ell) L_{\cG} + O(h)$.
The hidden constants depend only on $\ell, L_{\cG}, c_2, h_0$ and $t_n$ but not on $h, n$ nor the curvature of the manifold $\Mr$.
\end{theorem}

\begin{proof}
Let us omit the notation $\cdot^{\mathrm{R}}$. As in the proof of Theorem \ref{theorem: projected exponential Euler convergence}, we start by applying Lemma \ref{lemma:retraction inequality}:
\begin{equation} \label{eq: proj exp runge retraction inequality}
\norm{Y_{n+1} - X(t_{n+1})} \leq \norm{\Tr(\phi_{P \cF}^h (Y_n) + K_{n+1}) - \phi_{P \cF}^h (Y_n) + M_{n+1}} \leq 2 \norm{K_{n+1}} + \norm{M_{n+1}},
\end{equation}
where, with the help of Lemma \ref{lemma: expansion of the projected flow}, we have
$$M_{n+1} = \phi_{\cF}^h (X(t_n)) - \phi_{P \cF}^h (Y_n),$$
$$K_{n+1} = h \sum_{i=1}^2 b_i(h \cL) \left[\proj{Y_{ni}}{\cG(Y_{ni})} - f(t_n + c_i h) \right] - \delta_{n+1} - \omega_{n+1} - \lambda_{n+1}.$$
The first difference term is, again, easily bounded by Theorem \ref{thm:DLRA theoretical error}:
\begin{equation} \label{eq: proj exp runge bound M}
\norm{M_{n+1}} \leq e^{h \ellF} \norm{X(t_n) - Y_n} + h \varphi_1(h \ell_F) \cdot \varepsilon_r.
\end{equation}
Bounding the second term involves the stage errors $E_{ni}^Y = Y_{ni} - X(t_n + c_i h)$ for $i = 1, 2$. From the definition of the stages, we immediately get that $\norm{E_{n1}^Y} = \norm{Y_n - X(t_n)}$ and
\begin{equation}
\norm{E_{n2}^Y} = \norm{Y_{n2} - X(t_n+c_2 h)} = \norm{\Tr(\phi_{P \cF}^{c_2 h} (Y_n) + K_{n2}) - \phi_{P \cF}^{c_2 h} (Y_n) + M_{n2}} \leq 2 \norm{K_{n2}} + \norm{M_{n2}},
\end{equation}
where $K_{n2}$ and $M_{n2}$ are difference terms.
Here, we could use the expansion on the stages derived and bounded in Lemma \ref{lemma: expansion of the projected flow}.
More concise is to notice that $Y_{n2}$ coincides with exponential Euler with a step size of $c_2 h$, so its error bound is already derived in the course of the proof of Theorem \ref{theorem: projected exponential Euler convergence}:
\begin{equation} \label{eq: proj exp runge bound En2}
\begin{aligned}
\norm{E_{n2}^Y} &\leq \left( e^{c_2 h \ellF} + c_2 h \varphi_1(c_2 h \ell) L_{\cG} \cdot (4 + O(c_2 h)) \right) \norm{E_n^Y} + (c_2 h)^2 \varphi_2(c_2 h \ell) \cdot \max_{0 \leq s \leq c_2 h} \norm{f'(t_n + s)} + C \cdot h \cdot \varepsilon_r \\
&\leq (1 + c_2 h \tilde{L}) \norm{E_n^Y} + (c_2 h)^2 \varphi_2(c_2 h \ell) \cdot \max_{0 \leq s \leq c_2 h} \norm{f'(t_n + s)} + C \cdot h \cdot \varepsilon_r,
\end{aligned}
\end{equation}
where $\tilde{L} = \ellF + 4 \varphi_1(c_2 h \ell) L_{\cG} + O(c_2 h)$.
By the triangular inequality, we have that
\begin{equation} \label{eq: initial bound Kn}
\norm{K_{n+1}} \leq h \sum_{i=1}^2 \norm{b_i (h \cL)} \cdot \left( L_{\cG} \norm{Y_{ni} - X(t_n + c_i h)} + \varepsilon_r \right) + \norm{\delta_{n+1}} + \norm{\omega_{n+1}} + \norm{\lambda_{n+1}},
\end{equation}
where we have the following bounds:
\begin{equation} \label{eq: proj exp runge bound delta omega lambda}
\begin{split}
\norm{\delta_{n+1}} \leq \tilde C_2 \cdot h^3 \cdot \max_{0 \leq s \leq h} \norm{f''(t_n + s)}, \quad
\norm{\omega_{n+1}} \leq h \varphi_1(h \ell) \varepsilon_r, \\
\text{and } \quad
\norm{\lambda_{n+1}} \leq L_{\cG} \cdot h \varphi_1(h \ell) \cdot (1 + O(h)) \cdot \norm{E_n} + O(h^2) \cdot \varepsilon_r.
\end{split}
\end{equation}
The end of the proof consists essentially in assembling the terms together.
Recall the bounds for $b_1, b_2$ from \eqref{eq: bounds for bi} and plugging \eqref{eq: proj exp runge bound En2} into \eqref{eq: initial bound Kn} gives, after simplification, the following bound:
\begin{equation} \label{eq: proj exp runge bound Kn}
\norm{K_{n+1}} \leq h \cdot L_{\cG} \cdot \left( 2 \varphi_1(h \ell) + O(h) \right) \norm{E_n^Y} + \tilde C_1 \cdot h^3 \left( \max_{0 \leq s \leq h} \norm{f'(t_n+s)} + \max_{0 \leq s \leq h} \norm{f''(t_n+s)} \right) + \tilde C_2 \cdot h \varepsilon_r.
\end{equation}
Now inserting \eqref{eq: proj exp runge bound Kn} and \eqref{eq: proj exp runge bound M} into \eqref{eq: proj exp runge retraction inequality} gives the following recursion for the error,
$$\norm{Y_{n+1} - X(t_{n+1})} \leq (1 + h L_*) \norm{Y_n - X(t_n)} + h (\beta + \gamma),$$
where we have defined the quantities $L_* = \ellF + 4 \varphi_1(h \ell) L_{\cG} + O(h),$ and
$$\beta = \hat C_1 \cdot h^2 \cdot  \left( \max_{0 \leq s \leq t_{n+1}} \norm{f'(s)} + \max_{0 \leq s \leq t_{n+1}} \norm{f''(s)} \right), \quad
\gamma = \hat C_2 \cdot \varepsilon_r.$$
The conclusion follows since solving the recursion leads to
$$\norm{Y_n - X(t_n)} \leq e^{nh L_*} \norm{Y_0 - X_0} + \frac{e^{nh L_*} - 1}{L_*} ( \beta + \gamma),$$
and keeping track of the hidden constants shows that they only depend on $\ell, L_{\cG}, c_2$ and $h_0$.
\end{proof}

\section{Efficient implementation with Krylov subspace techniques} \label{sec: efficient implementation}

Having proved the global error of the projected exponential Euler and Runge schemes, we will now explain their efficient implementation. As usual, the presence of the $\varphi$-functions requires special attention for large-scale problems. An added complication in our setting is that all computations should be done with low-rank approximations in mind. For example, storing or computing full matrices has to be avoided. Most of this section is therefore devoted in explaining how Krylov subspace approximation techniques for $\varphi$-functions can be performed efficiently with low-rank matrix approximations.

The efficient computation of $\varphi_k(hM)x$ given a vector $x$ and a matrix $M$ by Krylov subspace techniques is a well-studied topic. In case of the matrix exponential $\varphi_0$, the earliest results go back to~\cite{gallopoulos1992efficient, hochbruck1997krylov}. The approximation of general $\varphi$-functions by Krylov spaces is a more recent topic discussed in \cite{niesen2009krylov} and by rational Krylov spaces in \cite{gockler2014uniform,gockler2014rational}. Adaptive techniques are also described in \cite{gaudreault2018kiops,botchev2020art}. 
Good software implementations are also available. For standard Krylov subspaces, Expokit \cite{sidje1998expokit} allows computing the action of the matrix exponential. It is originally written in FORTRAN and \textsc{Matlab}, and has been embedded into other languages, including Python. For rational Krylov methods, there exists the RKToolbox \cite{beckermann2010convergence,elsworth2020block,ruhe1994rational,ruhe1998rational,wilkinson1965algebraic}, which is written in \textsc{Matlab}. Finally, when the matrix $M$ is relatively small, dense linear algebra techniques for the matrix exponential are preferable and we refer to \cite{al2011computing} for an efficient implementation in the context of  $\varphi$-functions.

The techniques mentioned above can be used directly to implement the projected exponential methods proposed in this paper. However, a good implementation should exploit that the $\varphi$-functions are evaluated for $\cL = I \otimes A + B^T \otimes I$, which has an important Kronecker structure. This is crucial when we want to obtain efficient low-rank approximations. 
The main idea is to build two independent Krylov subspaces for the range and co-range of the low-rank approximation. This is similar to Krylov methods that compute a low-rank approximation for \emph{algebraic} Lyapunov and Sylvester equations; see, e.g., \cite{simoncini2007new} for extended Krylov subspaces and \cite{simoncini2016computational} for a recent overview. Applied to differential Sylvester equations (hence, the ODE~\eqref{eq:full problem} with a constant function $\cG$), such standard Krylov subspaces were already used in \cite{behr2019solution}. In \cite{koskela2017analysis}, this was generalized to differential Riccati equations (hence, with $\cG(X) = Q - X(t)SX(t)$) and analyzed for standard Krylov subspaces. The techniques proposed here find similarities with these methods, but our approach treats the differential equation~\eqref{eq:full problem} with general $\cG$ and high-order exponential integrators. In addition, we analyze the error for rational Krylov subspaces.

\subsection{Rational Krylov approximation} 
Before going into the details of the new technique, let us introduce basic notions of Krylov spaces. Given a generic (sparse) matrix $S \in \Rmn$ and a (tall) matrix $X \in \Rnr$, we define three variants of (block) Krylov spaces:
\begin{align*}
     & \text{Polynomial Krylov space:}          &  & \mathcal K_k(S, X) = \mathrm{span} \{X, SX, S^2X, \ldots, S^{k-1}X \},                                                                                 \\
     & \text{Extended Krylov space:} &  &E \mathcal{K}_k(S, X) = \mathcal K_k(S, X) + \mathcal K_k(S^{-1}, S^{-1} X),                                                                            \\
     & \text{Rational Krylov space:} &  &R \mathcal{K}_k(S, X) = q_{k-1}(S)^{-1} \mathcal K_k(S, X), \quad \text{for a degree $k-1$ polynomial $q_{k-1}$ }.
\end{align*}
An extended Krylov space requires that $S$ is invertible. For the construction of the rational Krylov space, one typically uses the factored polynomial $q_{k-1}(z) = (z-\rho_1) \cdots (z-\rho_{k-1})$. To ensure the invertibility of $q_{k-1}(S)$, the poles $\rho_1, \ldots, \rho_{k-1}$ have to be different from the eigenvalues of $S$. Observe that the poles are part of the definition of $R \mathcal{K}_k(S, X)$ and can be fixed from the start or adaptively during the iteration. We will focus on the former for simplicity.

Computing orthonormal bases for these block Krylov spaces requires some care but can be done efficiently with the (rational) Arnoldi algorithm; we refer to \cite{berljafa2014rational} for a \textsc{Matlab} implementation. Good (and even optimal) choices of the poles are discussed in \cite{guttel2013rational}.

In the following, a generic Krylov space $G \mathcal K_k(s,X)$ will refer to any of these spaces.
The matrix $S$ represents a generic matrix (operator) but will always be $A$ or $B$ from~\eqref{eq:full problem} in the following.
Polynomial Krylov spaces are the least expensive to compute per step but they also have the worst approximation error leading to many iterations and thus higher memory requirements. Rational Krylov spaces have the best approximation error but they require determining the poles and their application per step is expensive. On the other hand, they typically require less iterations and, therefore, have a low memory footprint. Extended Krylov spaces are the compromise between these two techniques if linear systems with $S$ are feasible to compute.

While any of these methods can be used, we will focus our analysis on rational Krylov spaces since they are the most memory efficient. 
The following lemma will be useful later in the analysis. It relates the rational Krylov approximation of $e^{tS}X$ to that of the best rational approximant of $e^{t}$. The latter has a long history and many results are known depending on the way the poles are determined; see the references in \cite{guttel2013rational}. We have chosen to give one of the more elegant results when $S$ is a  symmetric matrix with strictly negative eigenvalue, hence $\ell_S < 0$. In the context of stiffness, this situation is actually very relevant. One can allow non-symmetric matrices as well using the field of values, but we do not give details for simplicity.

\begin{lemma}[Error made by rational Krylov approximation on the exponential function] \label{lemma:error made by rational Krylov exponential} 
    Consider a symmetric matrix $S \in \Rnn$ with $\ell_S = \lambda_{\max}(S) < 0$ and a matrix $X \in \mathbb{R}^{n \times m}$. Let $Q_k$ be such that $Q_k^T Q_k = I$ and $\myspan(Q_k) = R \mathcal{K}_k(S, X)$ with $R\mathcal{K}_k(S, X)$ the rational Krylov space with optimal poles $\rho_1, \ldots, \rho_{k-1}$ as defined in the proof. Denote $S_k = Q_k^T S Q_k$. Then, for all $t \geq 0$ it holds
\[
  \norm{e^{t S}X - Q_k e^{t S_k} Q_k^T X} \leq \frac{2 \sqrt{2}}{3^{k-1}} \norm{X}.
\]
In the asymptotic regime $k\to \infty$, the constant $3$ can be improved to at least $9.037$.
\end{lemma}
\begin{proof}
Let $\Pi_{k-1}$ be the set of polynomials of degree at most $k-1$. By Lemma 3.2 in \cite{beckermann2021low}, the block rational Krylov approximation has the exactness property
\[
 r(S)X = Q_k \, r(Q_k^T S Q_k) \, Q_k^T X \qquad \forall r = p/q_{k-1}, \ p \in \Pi_{k-1}.
\]
Here, $q_{k-1} \in \Pi_{k-1}$ depends on the poles that we leave undetermined for now. From this and the fact that $R\mathcal{K}_k(S, X)=R\mathcal{K}_k(tS, X)$ for nonzero $t$, we obtain that the rational Krylov approximation error of the exponential is never worse than that by any rational function $r$ as defined above,
\begin{align*}
 \| e^{t S}X - Q_k e^{t S_k} Q_k^T X \| &= \| e^{t S}X - r(tS) X + Q_k \, r(tQ_k^T S Q_k) Q_k^T X - Q_k e^{t S_k} Q_k^T X \| \\
 &\leq  \left( \| e^{t S} - r(tS) \|_2 + \| e^{t S_k} - r(tQ_k^T S Q_k) \|_2 \right) \|X\|.
\end{align*}
Since $S$ and $Q_k^T S Q_k$ are symmetric, the approximation errors above can be calculated on their spectrum. Both matrices having negative eigenvalues and $t \geq 0$, we can therefore bound
\begin{equation}\label{eq:approx Krylov to scalar}
 \| e^{t S}X - Q_k e^{t S_k} Q_k^T X \| \leq 2 \, \|X\| \, \sup_{z \in [0, \infty)} | e^{-z} - r(-z) |.
\end{equation}
Recall that $r = p/q_{k-1}$ such that $p \in \Pi_{k-1}$ was arbitrary and $q_{k-1} \in \Pi_{k-1}$ depends on the still to be fixed poles. Hence, we can take the infimum of~\eqref{eq:approx Krylov to scalar} over all such rational functions $r$. By classical results  in approximation theory (see, e.g., \cite[Thm.~4.11]{petrushev2011rational}), we have
\[
 \inf_{r}\sup_{z \in [0, \infty)} | e^{-z} - r(-z) | = \frac{\sqrt{2}}{3^{k-1}}.
\]
Combining leads to the desired result with the poles taken from the best approximant $r_*$ above. The statement about the asymptotic improvement is also in
\cite[Chap.~4.5]{petrushev2011rational} and \cite{varga1993numerical}.
\end{proof}

\subsection{Lucky Krylov approximation}

In this section, we explain how to exploit the structure of the projected exponential methods \eqref{eq: projected exponential methods} and perform efficient computations.
In \cite[Lemma 4.1]{koch2007dynamical}, it is shown that the projection satisfies
\begin{align} \label{eq: projection formula}
    \proj{Y}{\cG(Y)} = UU^T \cG(Y) - UU^T \cG(Y) VV^T + \cG(Y) VV^T,
\end{align}
where $Y = U \Sigma V^T$ is a compact SVD decomposition. The key idea in our implementation is to exploit that $Y$ can be represented in the same subspace that is used for $\proj{Y}{\cG(Y)}$, independent of $\mathcal G$. Together with Lemma \ref{lemma: phi functions as IVPs}, we will show that the sub-steps of the projected exponential methods can then be computed all-at-once in an efficient manner with Krylov techniques.

\subsubsection*{First order: projected exponential Euler}

Recall the projected exponential Euler scheme \eqref{eq:projected exponential euler}:
\begin{align}\label{eq:exp Euler with temp Y1}
    Y_1 = \Ret{\tilde{Y}_1} = \Ret{ e^{h \cL} Y_0 + h \varphi_1 (h \cL) \proj{Y_0}{\cG (Y_0)} }.
\end{align}
The naive approach would be to compute its vectorized formulation,
$$\mathrm{vec}(\tilde{Y}_1) = e^{h \cL} \mathrm{vec}(Y_0) + h \varphi_1 (h \cL) \mathrm{vec}(\proj{Y_0}{\cG (Y_0)}),$$
using classical Krylov techniques for exponential functions, as described in \cite{sidje1998expokit} and \cite{al2011computing}.
However, one would have to form dense matrices, which makes such an approach very inefficient.
Fortunately, it is possible to reduce the computational cost by smartly exploiting the Sylvester-like structure.

By Lemma \ref{lemma: phi functions as IVPs}, the term inside the truncation can be obtained as $Z(h)$ where $Z$ satisfies the following Sylvester differential equation:
\begin{align} \label{eq: equivalence for projected Euler}
    Z(t) = e^{t \cL} Y_0 + t \varphi_1 (t \cL) \proj{Y_0}{\cG (Y_0)} \iff
    \left\{
    \begin{aligned}
         & \dt{Z}(t) = A Z(t) + Z(t) B + \proj{Y_0}{\cG (Y_0)}, \\
         & Z(0) = Y_0.
    \end{aligned}
    \right.
\end{align}
Since $Y_0 \in \Mr$ and $\proj{Y_0}{\cG (Y_0)} \in \mathcal T_{Y_0} \mathcal M_r$, these matrices can be represented using the SVD as
\begin{align*}
     & Y_0 = U_0 \ \Sigma \ V_0^T                                      &  & \text{ where } \Sigma \in \Rrr, U_0 \in \Rnr, V_0 \in \Rmr;                      \\
     & \proj{Y_0}{G(Y_0)} = [U_0, U_1] \ \tilde{\Sigma} \ [V_0, V_1]^T &  & \text{ where } \tilde{\Sigma} \in \R^{2r \times 2r}, U_1 \in \Rnr, V_1 \in \Rmr.
\end{align*}
Observe that $Y_0$ can be represented in the bases for the range and co-range of $\proj{Y_0}{G(Y_0)}$. It is therefore natural to consider the following two Krylov spaces
\begin{equation*}
    \begin{aligned}
         & Q_k \text{ is a matrix with orthonormal columns s.t.}  &  \myspan(Q_k) = G\mathcal{K}_k (A, \mathrm{span}([U_0, U_1])), \\
         & W_k \text{ is a matrix with orthonormal columns s.t.}  &  \myspan(W_k) = G\mathcal{K}_k (B, \mathrm{span}([V_0, V_1])), \\
    \end{aligned}
\end{equation*}
Then, we perform a Galerkin projection to obtain a reduced Sylvester differential equation,
\begin{align} \label{eq: reduced Sylvester IVP}
    \left\{
    \begin{aligned}
         & \dt{S_k}(t) = Q_k^T A Q_k S_k(t) + S_k(t) W_k^T B W_k + Q_k^T \proj{Y_0}{\cG (Y_0)} W_k \\
         & S_k(0) = Q_k^T Y_0 W_k
    \end{aligned}
    \right.
\end{align}
where $S_k(t) \in \R^{2kr \times 2kr}$.
Let us denote $A_k = Q_k^T A Q_k$ and $B_k = W_k^T B W_k$.
When the Sylvester operator is invertible, the closed-form solution is given by (see, e.g., \cite{behr2019solution})
\begin{equation*}
S_k(t) = e^{t A_k} (S_k(0) + C) e^{t B_k} - C, \qquad  A_k C + C B_k = Q_k^T \proj{Y_0}{\cG (Y_0)} W_k,
\end{equation*}
where the second equation is an algebraic Sylvester equation.
These computations are feasible since the dimension is much smaller than in \eqref{eq: equivalence for projected Euler}. In particular, it can be evaluated in $O(k^3r^3)$ flops with techniques from dense linear algebra for computing the matrix exponential~\cite{al2010new} and solving the algebraic Sylvester equation~\cite{bartels1972solution}. Alternatively, any integrator for stiff ODEs can be used for computing the solution to the reduced IVP.
Finally, with the approximated solution $S_k(h)$ of \eqref{eq: equivalence for projected Euler} we can define 
\begin{equation}\label{eq:exp Euler approx Krylov Zh}
\tilde Y_1 = Z(h) \approx Z_k(h) = Q_k S_k(h) W_k^T.
\end{equation}

\subsubsection*{Second order: projected exponential Runge}

Let us recall the second-order scheme \eqref{eq: proj exp Runge} with $c_2 = 1$:
\begin{equation*} 
\begin{aligned}
    Y_{1/2}^{\mathrm{PR}} & = \Ret{e^{h \cL} Y_0 + h \varphi_1 \left( h \cL \right) \proj{Y_0}{\cG(Y_0)}}, \\
    Y_1^{\mathrm{PR}}     & = \Ret{e^{h \cL} Y_0 + h \varphi_1 \left( h \cL \right) \proj{Y_0}{\cG(Y_0)} + h \varphi_2(h \cL) \left( \proj{Y_{1/2}^{\mathrm{PR}}}{\cG(Y_{1/2}^{\mathrm{PR}})} - \proj{Y_0}{\cG(Y_0)} \right) }.
\end{aligned}
\end{equation*}
The term $Y_{1/2}^{\mathrm{PR}}$ can be approximated with the first order approximation described above. 
Similarly, we use Lemma \ref{lemma: phi functions as IVPs} for treating the term $Y_1^{\mathrm{PR}}$. The term inside the truncation is equivalent to $Z(h)$ where $Z$ satisfies the initial value problem
\begin{equation} \label{eq: system order 2}
\left\{
\begin{aligned}
&\dt{Z}(t) = A Z(t) + Z(t) B + \proj{Y_0}{\cG(Y_0)} + \frac{t}{h} \left[ \proj{\Yrungehalf}{\cG(\Yrungehalf)} - \proj{Y_0}{\cG(Y_0)} \right], \\
&Z(0) = Y_0.
\end{aligned}
\right.
\end{equation}
By the SVD, we can write
\begin{equation*}
Y_0 = U_0 \Sigma V_0^T, \quad \proj{Y_0}{\cG(Y_0)} = [U_0, U_1] \tilde{\Sigma} [V_0, V_1]^T, \quad \proj{\Yrungehalf}{\cG(\Yrungehalf)} = U_2 \hat{\Sigma} V_2^T,
\end{equation*}
for some matrices $\Sigma \in \Rrr, \tilde{\Sigma} \in \R^{2r \times 2r}$, and $\hat{\Sigma} \in \R^{2r \times 2r}$.
Clearly all three matrices can be represented in the subspaces $[U_0, U_1, U_2]$ and $[V_0, V_1, U_2]$, each of dimension $4k$. While in general $\myspan(U_2)~\not\subset~\myspan([U_0, U_1])$, these subspaces are close for small step size $h$. It is therefore reasonable to truncate these approximation subspaces to lower dimension in practice. In any case  the structure suggests using two Krylov spaces such that
\begin{equation*}
    \begin{aligned}
         & Q_k \text{ is a matrix with orthonormal columns s.t.}  &  &\myspan(Q_k) = G\mathcal{K}_k (A, \myspan([U_0, U_1, U_2])), \\
         & W_k \text{ is a matrix with orthonormal columns s.t.} &  &\myspan(W_k) = G\mathcal{K}_k (B, \myspan([V_0, V_1, U_2])).
    \end{aligned}
\end{equation*}
The reduced system is then obtained by a Galerkin projection,
\begin{equation} \label{eq: reduced system order 2}
\left\{
\begin{aligned}
&\dt{S}_k(t) = A_k S_k(t) + S_k(t) B_k + Q_k^T \proj{Y_0}{\cG(Y_0)} W_k + \frac{t}{h} Q_k^T \left[ \proj{\Yrungehalf}{\cG(\Yrungehalf)} - \proj{Y_0}{\cG(Y_0)} \right] W_k, \\
&S_k(0) = Q_k^T Y_0 W_k,
\end{aligned}
\right.
\end{equation}
where $A_k = Q_k^T A Q_k, B_k = W_k^T B W_k$, and $S_k(t) \in \R^{4kr \times 4kr}$. Again, when the Sylvester operator is invertible, the closed form solution is given by
\begin{equation*}
\begin{aligned}
&h \cdot ( A_k \hat{D} + \hat{D} B_k ) = D, \quad A_k D + D B_k = Q_k^T \left[ \proj{\Yrungehalf}{\cG(\Yrungehalf)} - \proj{Y_0}{\cG(Y_0)} \right] W_k, \\
&A_k C + C B_k = Q_k^T \proj{Y_0}{\cG(Y_0)} W_k, \\
&S_k(t) = e^{t A_k} \left(S_k(0) + C + \hat{D} \right) e^{t B_k} - C - \hat{D} - D,
\end{aligned}
\end{equation*}
where the matrices $C, \, D$, and $\hat D$ are obtained by solving algebraic Sylvester equations at a cost of $O(k^3 r^3)$ flops. Alternatively, the reduced system \eqref{eq: reduced system order 2} can be solved by any numerical integrator for stiff ODEs. The final approximation is 
\begin{equation}\label{eq:exp Runge approx Krylov Zh}
\tilde Y_1^{\mathrm{PR}} = Z(h) \approx Z_k(h) = Q_k S_k(h) W_k^T.
\end{equation}

\begin{remark}
The procedure immediately extends to any higher-order scheme with $s$ stages. It requires computing $2s$ different Krylov spaces in total. The largest reduced system will have size $2skr$, at maximum leading to a cost of $O(s^3k^3r^3)$ flops. In practice, it is likely that the subspaces of the sub-steps are not strictly orthogonal leading to lower dimensional systems.
\end{remark}

\subsection{A priori approximation bounds}

The following results give a priori convergence bounds for the approximation error when rational Krylov spaces are used.
These bounds can be used for estimating the size of the Krylov spaces (that is, the number of iterations in the Arnoldi iteration). A similar bound is derived in \cite{koskela2017analysis} for polynomial Krylov spaces. Our results extend it to rational Krylov spaces and to higher-order schemes.

\begin{theorem}[Rational Krylov approximation error, Euler scheme] \label{thm: rational error euler}
    Let $Z(t)$ be the solution to the full differential Sylvester equation \eqref{eq: equivalence for projected Euler} with symmetric and negative definite matrices $A$ and $B$. Let $Z_k(t)~=~Q_k S_k(t) W_k^T$ with $S_k(t)$ be the solution to the reduced differential Sylvester equation \eqref{eq: reduced Sylvester IVP} obtained with a rational Krylov space with $k-1$ poles as in Lemma \ref{lemma:error made by rational Krylov exponential}.
    Then, the error of the approximation~\eqref{eq:exp Euler approx Krylov Zh} satisfies
    \begin{align*}
        \norm{Z(h) - Z_k(h)} \leq \frac{4 \sqrt{2}}{3^{k-1}} \left( e^{h \ell_*}  \norm{Y_0} + \frac{e^{h \ell_*} - 1}{\ell_*} \norm{P_{Y_0} \cG (Y_0)} \right) 
    \end{align*}
    where $\ell_* = \max (\ell_A, \ell_B) \leq 0$. Asymptotically, the constant $3$ improves to at least $9$. 
\end{theorem}

\begin{proof}
    The closed-form solutions satisfy
    \begin{align*}
        Z(t)   & = e^{tA} Y_0 e^{tB} + \int_0^t e^{(t-s)A} \proj{Y_0}{\cG(Y_0)} e^{(t-s) B} ds,                                                   \\
        Z_k(t) & = Q_k e^{t A_k} Q_k^T Y_0 W_k e^{t B_k} W_k^T + \int_0^t Q_k e^{(t-s)A_k} Q_k^T \proj{Y_0}{\cG(Y_0)} W_k e^{(t-s) B_k} W_k^T ds,
    \end{align*}
    where $A_k = Q_k^T A Q_k$ and $B_k = W_k^T B W_k$. Let us decompose the error as
    $$Z(t) - Z_k(t) = E_k^1(t) + E_k^2(t),$$
    where
    \begin{align*}
        E_k^1(t) & = e^{tA} Y_0 e^{tB} - Q_k e^{t A_k} Q_k^T Y_0 W_k e^{t B_k} W_k^T,                                                             \\
        E_k^2(t) & = \int_0^t e^{(t-s)A} \proj{Y_0}{\cG(Y_0)} e^{(t-s) B} - Q_k e^{(t-s)A_k} Q_k^T \proj{Y_0}{\cG(Y_0)} W_k e^{(t-s) B_k} W_k^T ds.
    \end{align*}
    Applying Lemma \ref{lemma:error made by rational Krylov exponential} twice, the norm of the first term is bounded as
    \begin{align*}
        \norm{E_k^1(t)} & = \norm{e^{tA} Y_0 e^{tB} - Q_k e^{t A_k} Q_k^T Y_0 e^{tB} + Q_k e^{t A_k} Q_k^T Y_0 e^{tB} - Q_k e^{t A_k} Q_k^T Y_0 W_k e^{t B_k} W_k^T} \\
                        & = \norm{\left( e^{tA} - Q_k e^{t A_k} Q_k^T \right) Y_0 e^{tB}} + \norm{Q_k e^{t A_k} Q_k^T Y_0 \left( e^{tB} - W_k e^{t B_k} W_k^T \right)} \\
                        & \leq  \frac{2 \sqrt{2}}{3^{k-1}} \cdot \|  Y_0 \| \cdot  \left(  \| e^{tB} \| + \| Q_k e^{t A_k} Q_k^T \| \right).
    \end{align*}
    With Lemma~\ref{lemma: phi functions are bounded} and $\ell_* = \max (\ell_A, \ell_B)$ we obtain
    \[
      \norm{E_k^1(t)} \leq  \frac{4 \sqrt{2}}{3^{k-1}} \cdot e^{t \ell_*} \norm{Y_0}.
    \]
    By a similar argument, the norm of the second term is bounded as follows
    \begin{align*}
        \norm{E_k^2(t)} &\leq \int_0^t  \frac{2 \sqrt{2}}{3^{k-1}} \cdot \| \proj{Y_0}{\cG(Y_0)} \| \cdot  \left(  \| e^{(t-s)B} \| + \| Q_k e^{(t-s) A_k} Q_k^T \| \right) ds \\
        & \leq \frac{4 \sqrt{2}}{3^{k-1}} \cdot \| \proj{Y_0}{\cG(Y_0)} \| \cdot \int_0^t e^{(t-s) \ell_*} ds \\
        & = \frac{4 \sqrt{2}}{3^{k-1}} \cdot \frac{e^{t \ell_*} - 1}{\ell_*}  \norm{\proj{Y_0}{\cG(Y_0)}}.
    \end{align*}
Substituting $t=h$ gives the desired result.
\end{proof}

Thanks to the fast exponential convergence in $k$, this theorem implies that the intermediate matrix $\tilde Y_1$ in the projected exponential Euler method~\eqref{eq:exp Euler with temp Y1} can be accurately computed in practice with only a few steps of the rational Krylov method. Since the rank of $\tilde Y_1$ is bounded by $2kr$, its rank will also not grow much and the truncation operator in~\eqref{eq:exp Euler with temp Y1} can be applied efficiently. Furthermore, the approximation error does not depend on the Lipschitz constant of $\mathcal L$ and is thus also robust to stiffness. 

The following result shows that essentially the same properties remain to hold for the projected exponential Runge scheme.

\begin{theorem}[Rational Krylov approximation error, Runge scheme] \label{thm: rational error runge}
Let $Z(t)$ be the solution to the full differential Sylvester equation \eqref{eq: system order 2} with symmetric and negative definite matrices $A$ and $B$. Let $Z_k(t)~=~Q_k S_k(t) W_k^T$ with $S_k(t)$ be the solution to the reduced differential Sylvester equation \eqref{eq: reduced system order 2} obtained with a rational Krylov space with $k-1$ poles as in Lemma \ref{lemma:error made by rational Krylov exponential}.
    Then, the error of the approximation~\eqref{eq:exp Runge approx Krylov Zh} satisfies
\begin{equation*}
\begin{split}
\norm{Z(h) - Z_k(h)} \leq \frac{4 \sqrt{2}}{3^{k-1}} \left( e^{h \ell_*}  \norm{Y_0} + \frac{e^{h \ell_*} - 1}{\ell_*} \norm{P_{Y_0} \cG (Y_0)} + \frac{e^{h \ell_*} - 1 - h \ell_*}{h \ell_*^2} \norm{D} \right),
\end{split}
\end{equation*}
where $\ell_* = \max (\ell_A, \ell_B) \leq 0$ and $D = \proj{\Yrungehalf}{\cG(\Yrungehalf)} - \proj{Y_0}{\cG(Y_0)}.$ Asymptotically, the constant $3$ improves to at least $9$. 
\end{theorem}

\begin{proof}
Let us define the two terms $E_k^1(t)$ and $E_k^2(t)$ as in the proof of Theorem~\ref{thm: rational error euler}. Writing down the closed-form formula for the solution $Z(t)$, we get the additional term
\begin{equation*}
E_k^3(t) = \frac{1}{h} \int_0^t s \left[ e^{(t-s)A} D e^{(t-s) B} - Q_k e^{(t-s)A_k} Q_k^T D W_k e^{(t-s) B_k} W_k^T \right] ds.
\end{equation*}
Taking its norm, we bound it like the terms $E_k^1(t)$ and $E_k^2(t)$ as follows
\begin{equation*}
\norm{E_k^3(t)} \leq \frac{4 \sqrt{2}}{3^{k-1}} \cdot \frac{1}{h} \int_0^t s e^{(t-s) \mu_*} \norm{D} ds = \frac{4 \sqrt{2}}{3^{k-1}} \frac{e^{t \mu_*} - 1 - t \mu_*}{h \mu_*^2} \norm{D}. \qedhere
\end{equation*} 
\end{proof}

\begin{remark}
The approximation bounds above extend easily to any linear combination of $\varphi$-functions and, therefore, to any higher-order (projected) exponential scheme.
\end{remark}

\section{Numerical experiments} \label{sec:numerical_experiments}

In this section, we test the proposed integrators numerically with stiff examples that contain a diffusion component. We start with the differential Lyapunov equation, where we show basic properties of the new methods. 
Then, we consider the Riccati differential equation, where we compare the new methods to existing methods. 
Moreover, we inspect the error made by the Krylov approximation techniques and compare them to the theoretical bounds derived in the previous section. 
Next, we demonstrate the capabilities of the new methods on the Allen-Cahn equation, when the non-linearity is not ``tangent space friendly''.
Finally, we propose a rank-adaptive variant of the methods and show its numerical behavior.

The implementation is done in Python and the experiments were performed on a MacBook Pro with M1 processor and 16GB of RAM, except for Figure \ref{fig:lyapunov errors} which required more memory and was performed on a workstation running Ubuntu 20.04 with an Intel i9-9900k processor and 64GB of RAM. All numerical experiments can be reproduced with the open source code available on \href{https://github.com/BenjaminCarrel/projected-exponential-methods}{GitHub}\footnote{The URL is: https://github.com/BenjaminCarrel/projected-exponential-methods}.

\subsection{Application to the differential Lyapunov equation}\label{sec:diff Lyap}

We consider the Lyapunov differential equation 
\begin{equation} \label{eq:lyapunov_diffeq}
    \dt{X}(t) = A X(t) + X(t) A^T + C(t), \quad X(0) = X_0,
\end{equation}
where $A \in \Rnn$ is sparse, and $C(t) \in \Rnn$ is a symmetric low-rank matrix.
The differential Lyapunov equation appears in discretized partial differential equations; but also in other areas, such as stability analysis, controller design for linear time-varying systems \cite{amato2014finite}, and optimal control of linear time-invariant systems on finite time horizons \cite{locatelli2002optimal}, to name a few. 
To demonstrate basic properties, we consider the partial differential equation modelling heat propagation\footnote{The heat propagation is the archetypal example of a stiff problem.} with Dirichlet boundary conditions,
\begin{equation} \label{eq:heat_propagation}
    \left\{ 
        \begin{aligned}
            \partial_t u(\pmb{x}, t) &= \Delta u(\pmb{x}, t) + s(\pmb{x}, t), &&\pmb{x} \in \Omega, \quad t \in [0, T], \\
            u(\pmb{x}, 0) &= \pmb{u}_0, &&\pmb{x} \in \Omega, \\
            u(\pmb{x}, t) &= 0, &&\pmb{x} \in \partial \Omega, \quad t \in [0, T].
        \end{aligned}
    \right.
\end{equation}
With two spatial dimensions, the discretized problem on a tensor grid (like standard finite differences) is exactly a Lyapunov differential equation where $A$ is the one-dimensional finite-difference stencil, and $C(t)$ is the discretization of the source $s(\pmb{x}, t)$. 
It is shown in \cite{carrel2023low} that, under a few conditions, the solution admits a good low-rank approximation. 

In the following, the problem is solved on the time-interval $[0, T] = [0, 1]$ and the spatial domain is $\Omega = [0, 1]^2$. 
If not mentioned otherwise, the size of the problem is $n=128$ (corresponding to $n^2 = 16384$ dofs), and the step size is $h=0.001$. 
For this mesh refinement, the problem is moderately stiff. 
An implicit method should be preferred, but an explicit method can still be used if the time steps are sufficiently small.
We use extended Krylov spaces with only one iteration, so $E \mathcal K(A, X) = \mathrm{span} \left\{X, A^{-1} X \right\}$.

The motivating example above, Figure \ref{fig:motivation}, was produced with a constant source term defined by the symmetric low-rank matrix
$$C(t) = Q \Sigma Q^T \in \R^{n \times n}, \quad \mathrm{diag}(\Sigma) = \left\{1, 10^{-4}, 10^{-8}, 10^{-12}, 10^{-16} \right\},$$
and $Q \in \R^{n \times 5}$ is a random matrix with orthonormal columns.
In this figure, we observe that the projected exponential Euler method is close to the approximation error. 
In fact, the source is constant so its derivative is zero and, therefore, Theorem \ref{theorem: projected exponential Euler convergence} predicts that the method is exact, up to the approximation error due to the rank and the Krylov space.

In the next experiments, we introduce a source term as defined in \cite{ostermann2019convergence}, but time-dependent. We consider the $q$ independent vectors $\left\{ \pmb{1}, \pmb{e}_1, \ldots, \pmb{e}_{(q-1)/2}, \pmb{f}_1, \ldots, \pmb{f}_{(q-1)/2} \right\},$ where 
$$e_k(x) = \sqrt{2} \cos (2 \pi k x) \quad \text{and} \quad f_k(x) = \sqrt{2} \sin (2 \pi k x), \quad k = 1, \ldots, (q-1)/2,$$
are evaluated at the grid points $\left\{ x_j \right\}_{j=1}^n$ with $x_j = \frac{j}{n+1}$. 
Then, the time-dependent source is defined as
\begin{equation*}
C(t) = e^{4 t} \cdot M^T M, \quad M = \begin{bmatrix} \pmb{1} & \pmb{e}_1 & \ldots & \pmb{e}_{(q-1)/2} & \pmb{f}_1 & \ldots & \pmb{f}_{(q-1)/2} 
\end{bmatrix}.
\end{equation*}
In the following, we take $q = 5$. The matrix $M$ will also be used in the Riccati experiments but with $q=9$.

\subsubsection*{Mesh refinements and performance}

\begin{figure}[!ht]
    \begin{subfigure}{0.49\textwidth}
        \includegraphics[width=\textwidth]{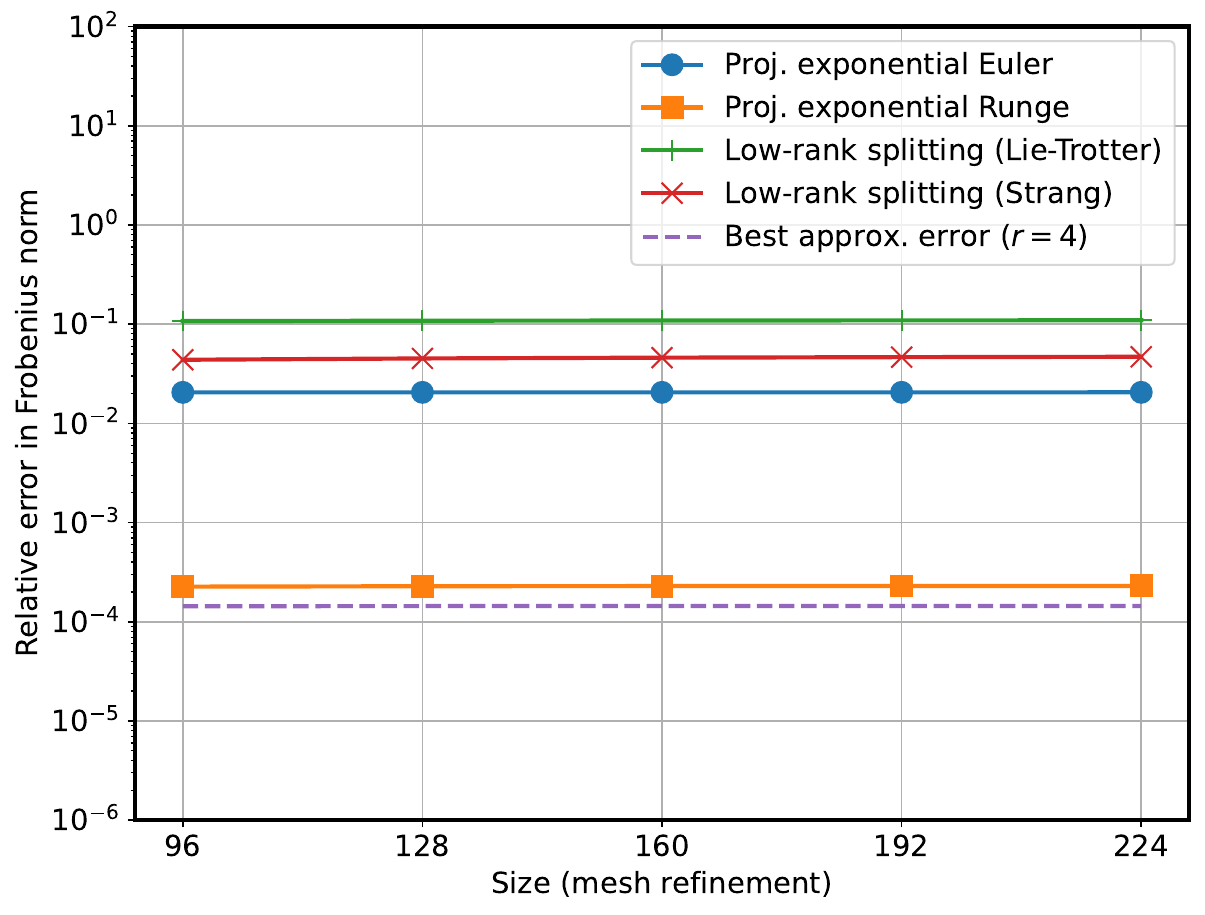}
        \subcaption{Relative error at final time under mesh refinements.}
    \end{subfigure}
      \hfill
    \begin{subfigure}{0.49\textwidth}
        \includegraphics[width=\textwidth]{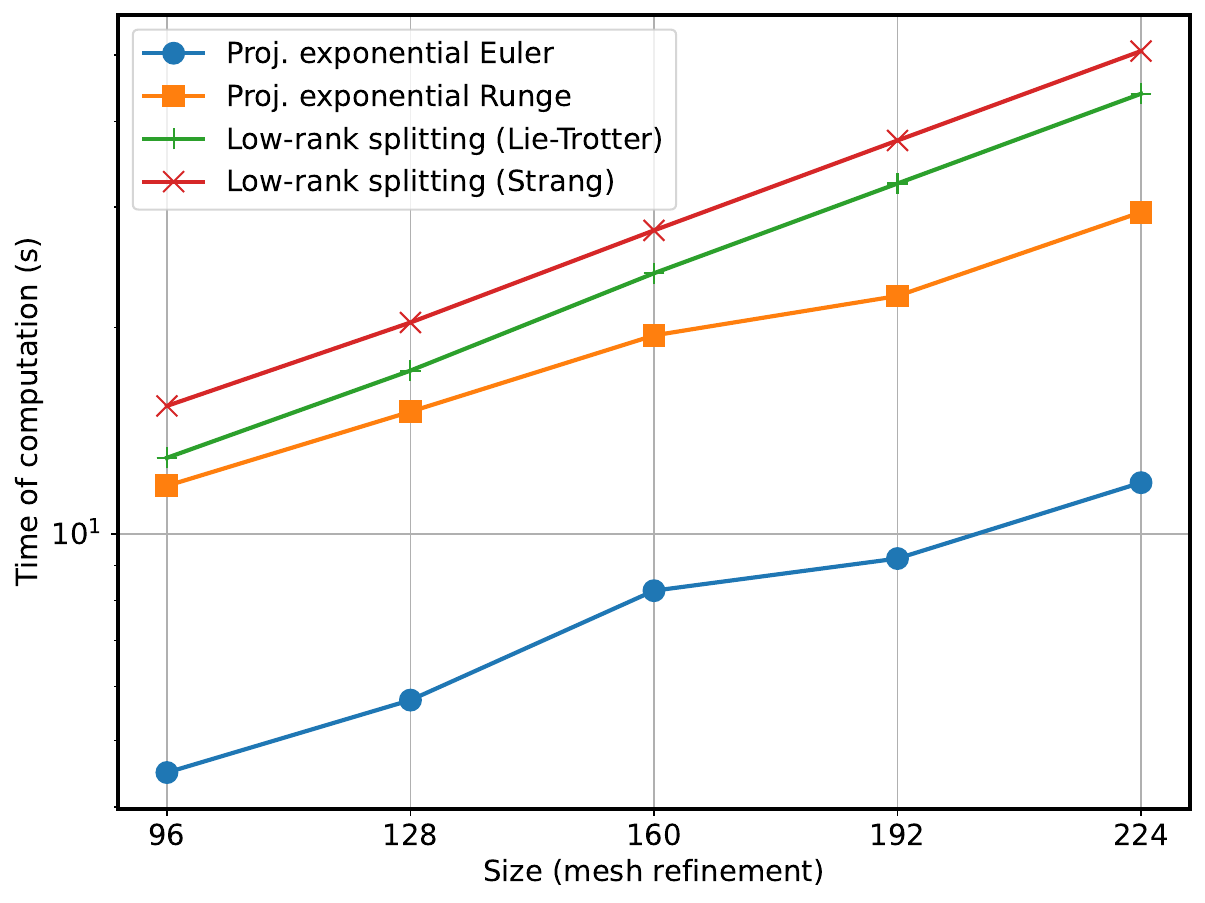}
        \subcaption{Total time of computation for each mesh size.}
    \end{subfigure}
    \centering
    \caption{New methods applied to the differential Lyapunov equation with time-dependent source. The low-rank splitting techniques are from \cite{ostermann2019convergence}.}
    \label{fig:lyapunov errors}
\end{figure}

In Figure \ref{fig:lyapunov errors}, we investigate the stability of the algorithms under mesh refinements. 
It is similar to the previous figure, but now with a time-dependent source and including second order methods. 
As we can see, our new methods are noticeably more accurate and also faster, in addition to being more precise than the existing techniques from \cite{ostermann2019convergence} based on low-rank splittings. 
Despite the mesh refinements; the errors remain stable, which verifies the robust-to-stiffness property of the methods.

\subsubsection*{Rank and performance}

\begin{figure}[!ht]
    \begin{subfigure}{0.49\textwidth}
        \includegraphics[width=\textwidth]{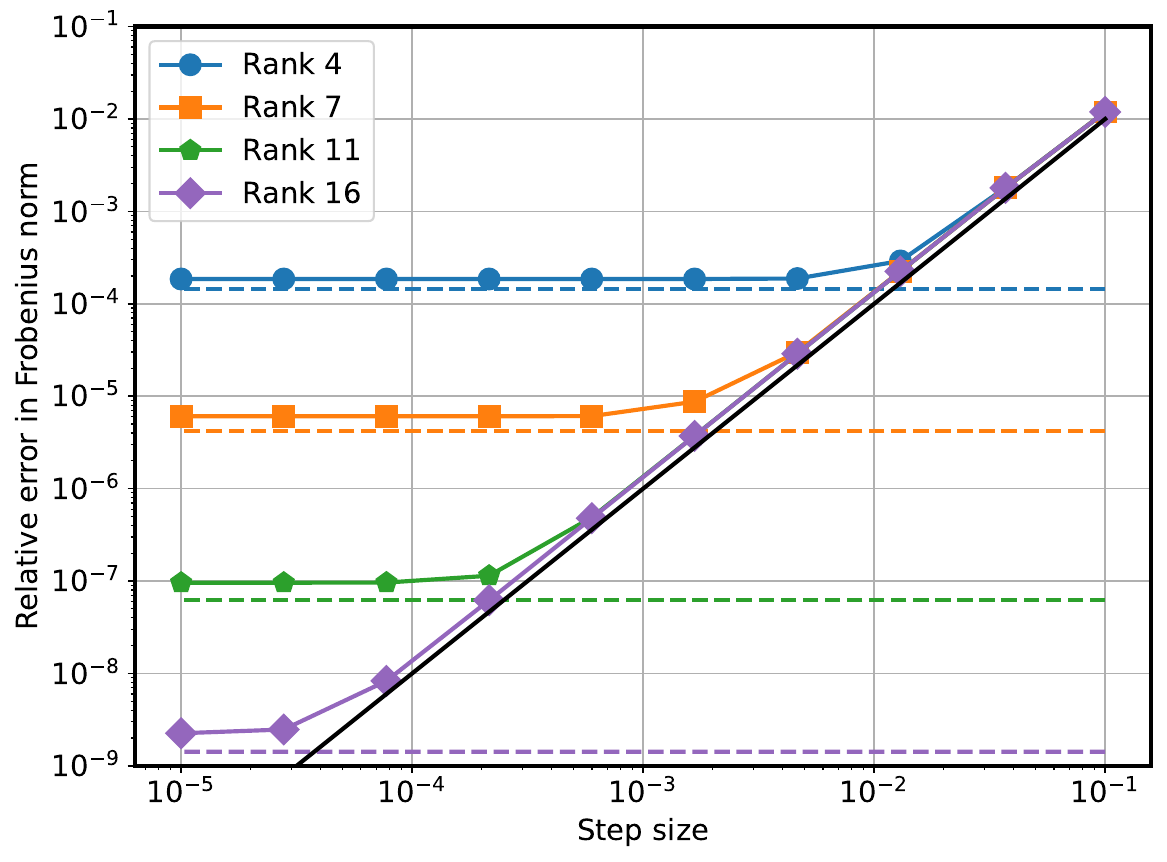}
        \subcaption{Relative error at final time for several ranks.}
    \end{subfigure}
      \hfill
    \begin{subfigure}{0.49\textwidth}
        \includegraphics[width=\textwidth]{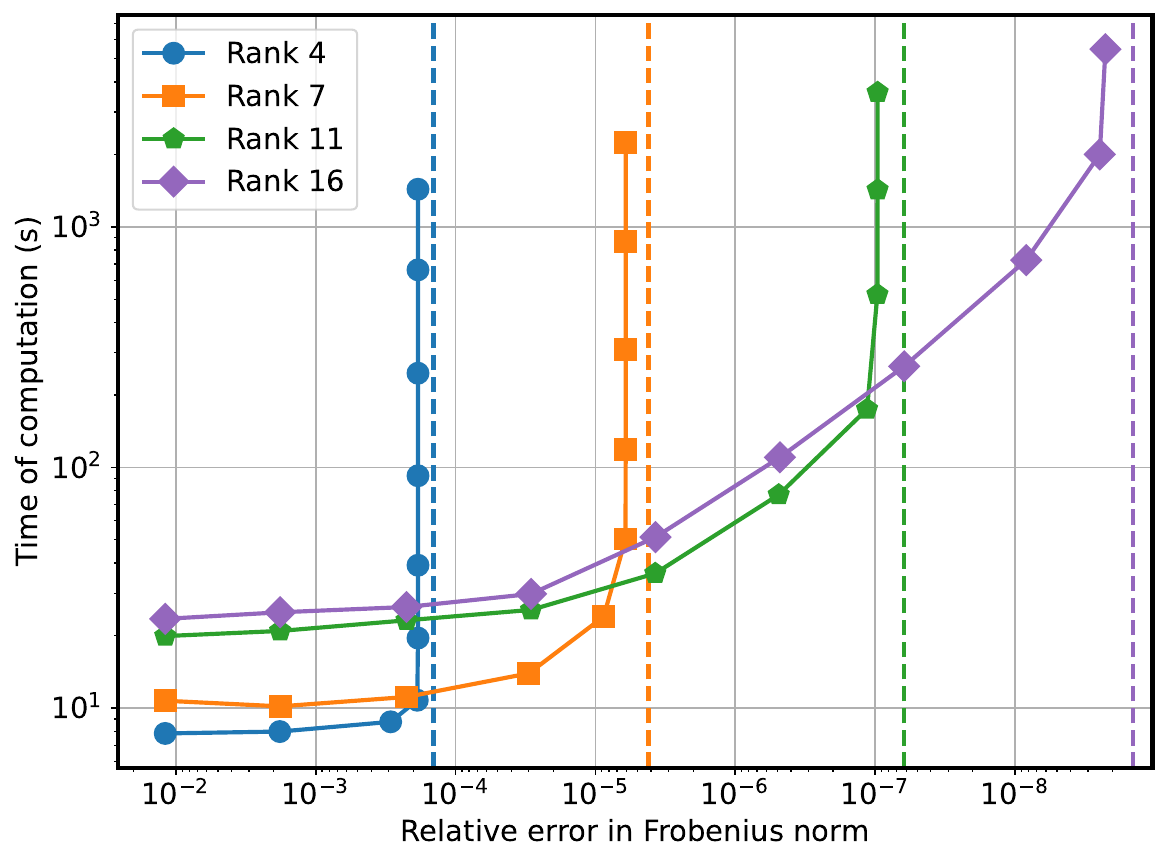}
        \subcaption{Total time of computation for each rank.}
    \end{subfigure}
    \centering
    \caption{Projected exponential Runge \eqref{eq: proj exp Runge} applied to differential Lyapunov equation with time-dependent source. The dashed lines indicate the minimal error for that particular rank.}
    \label{fig:lyapunov ranks}
\end{figure}

Another important property for a low-rank method is to be robust to small singular values.
In other words, the methods have to remain stable in the presence of small singular values, which is generally observed since we desire small error due to finite rank.
In Figure \ref{fig:lyapunov ranks}, we applied projected exponential Runge \eqref{eq: proj exp Runge} with several ranks of approximation.
Clearly, the convergence properties do not depend on the choice of the rank, which only limits the accuracy due to the low-rank approximation.
For this particular problem, the method shows very good performance, even when the rank increases.

\subsubsection*{Non strict order conditions}

\begin{figure}[!ht]
    \begin{subfigure}{0.49\textwidth}
        \includegraphics[width=\textwidth]{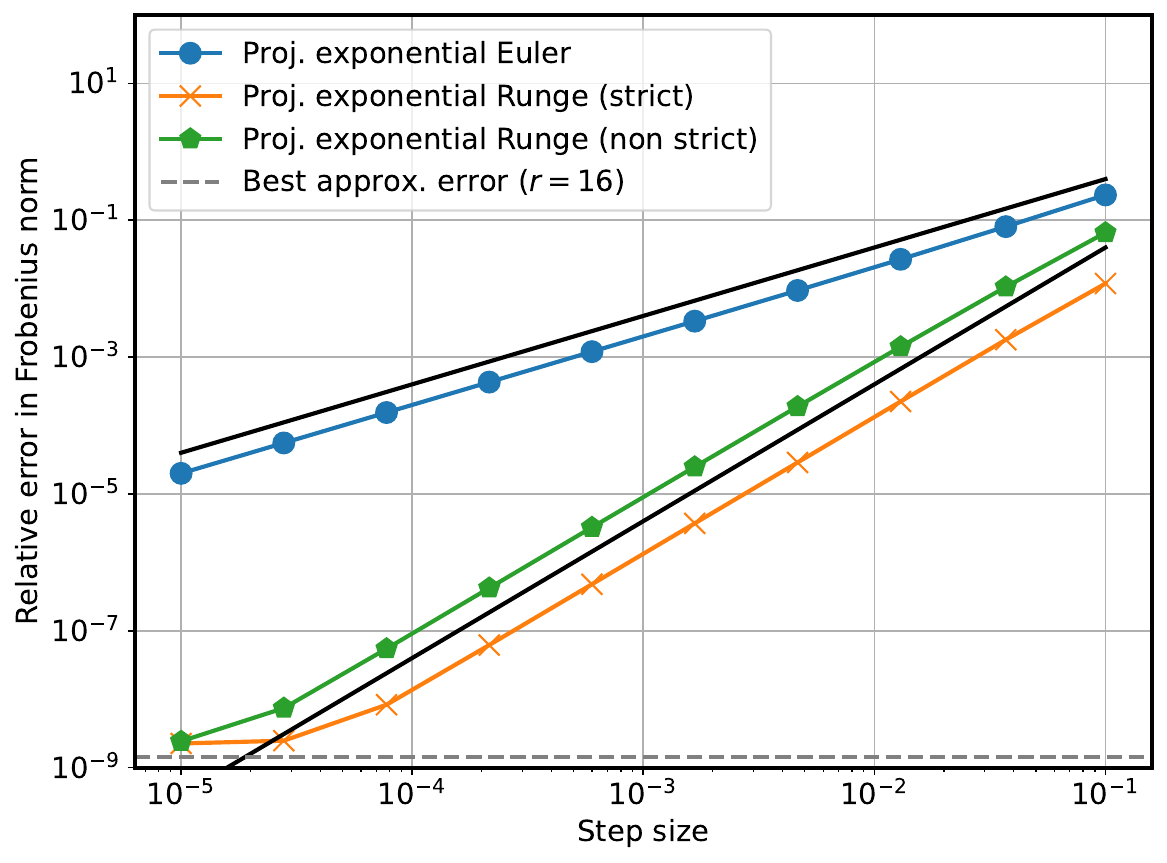}
        \subcaption{Relative error at final time.}
    \end{subfigure}
      \hfill
    \begin{subfigure}{0.49\textwidth}
        \includegraphics[width=\textwidth]{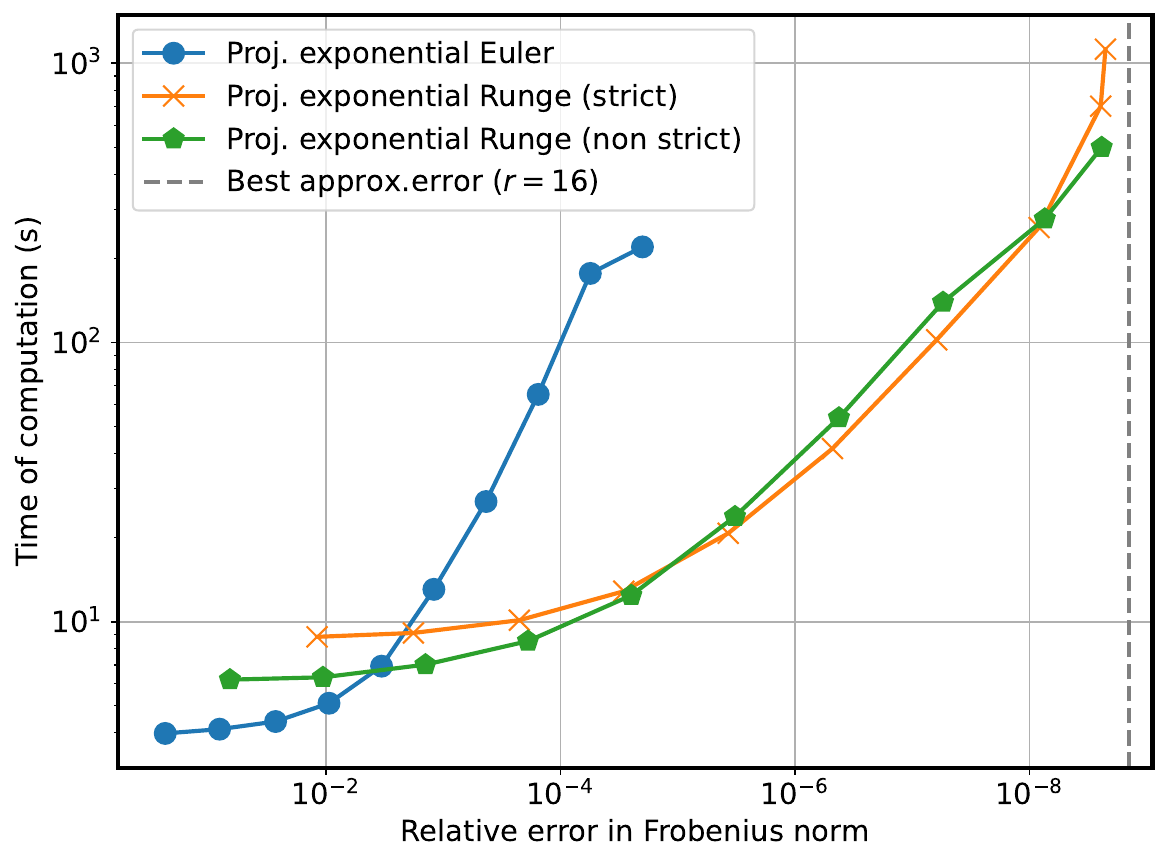}
        \subcaption{Total time of computation for each method.}
    \end{subfigure}
    \centering
    \caption{New methods applied to the differential Lyapunov equation with time-dependent source.}
    \label{fig:lyapunov three methods}
\end{figure}

As mentioned in Remark \ref{remark: non strict Runge}, weakening the second order condition leads to a method that has classical order two but not stiff order two. 
In Figure \ref{fig:lyapunov three methods}, we compare the error of projected exponential Runge with strict order conditions and its non strict equivalent. 
While the non-strict method does not have robust order two in theory, it performs very well in practice.
As we can see, the non-strict method has a slightly larger hidden constant than its strict equivalent, and the same convergence rate of two.
It is, therefore, a good practical compromise between the two other methods since, as mentioned earlier, it requires only evaluations of $\varphi_1$ functions. 
On the performance figure, it is clear that for a low accuracy, the projected exponential Euler method is faster. For moderate accuracy, one would prefer the non-strict projected exponential Runge; and we recommend using the strict projected exponential Runge for high accuracy.

\subsection{Application to the differential Riccati equation}

In order to compare with existing methods, we consider the problem described in \cite{ostermann2019convergence} which we briefly restate here.
The problem is a differential Riccati equation of the form
\begin{equation}\label{eq:riccati_diffeq}
  \dt{X}(t) = A^T X(t) + X(t) A + M^T M - X(t)^2, \quad X(0) = X_0.
\end{equation}
The sparse matrix $A \in \R^{n \times n}$ comes from the spatial discretization of the diffusion operator 
$$\mathcal D = \partial_x (\alpha(x) \partial_x(\cdot)) - \lambda I, \quad \alpha(x) = 2 + \cos(2 \pi x), \, \lambda = 1,$$ 
on the domain $\Omega = (0, 1)$ subject to homogeneous Dirichlet boundary conditions.
The discretization is done via the finite volume method, as described in \cite{gander2018numerical}. 
Similarly to before, the constant matrix $M \in \R^{q \times n}$ is defined by taking $q=9$ independent vectors $\left\{ 1, e_1, \ldots, e_{(q-1)/2}, f_1, \ldots, f_{(q-1)/2} \right\}.$
As initial value, we start from $\tilde{X}_0 = 0$, propagate the solution to time $\delta t=0.01$ with the reference solver in order to become a physically meaningful initial value $X_0$ that is full rank. 
The final time is $T = 0.1$ and the size is $n=200$.
As before, we use extended Krylov spaces with only one iteration.
The rank of approximation is $r=20$. 

\subsubsection*{Comparison with existing methods}

\begin{figure}[!ht]
    \begin{subfigure}{0.49\textwidth}
        \includegraphics[width=\textwidth]{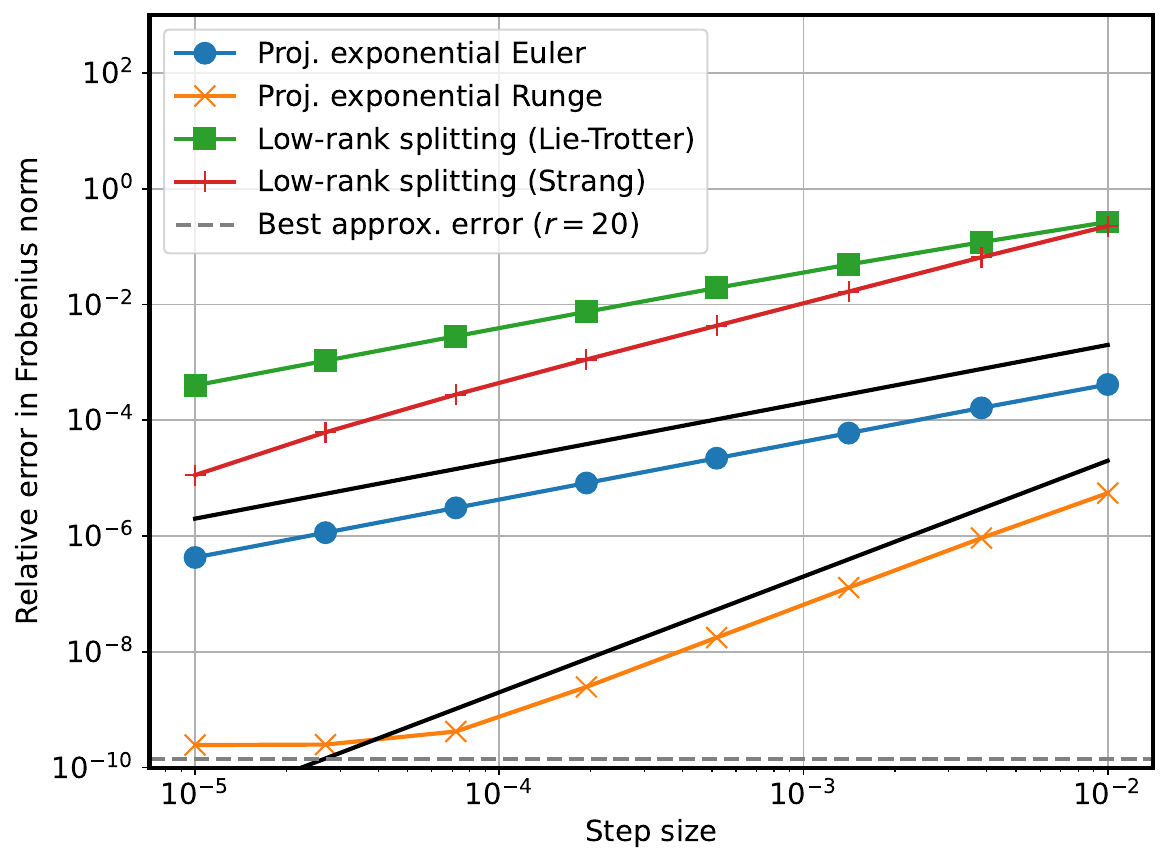}
        \subcaption{Global error. The straight black lines are slopes $1$ and $2$ for visual reference.}
    \end{subfigure}
    \hfill
    \begin{subfigure}{0.49\textwidth}
        \includegraphics[width=\textwidth]{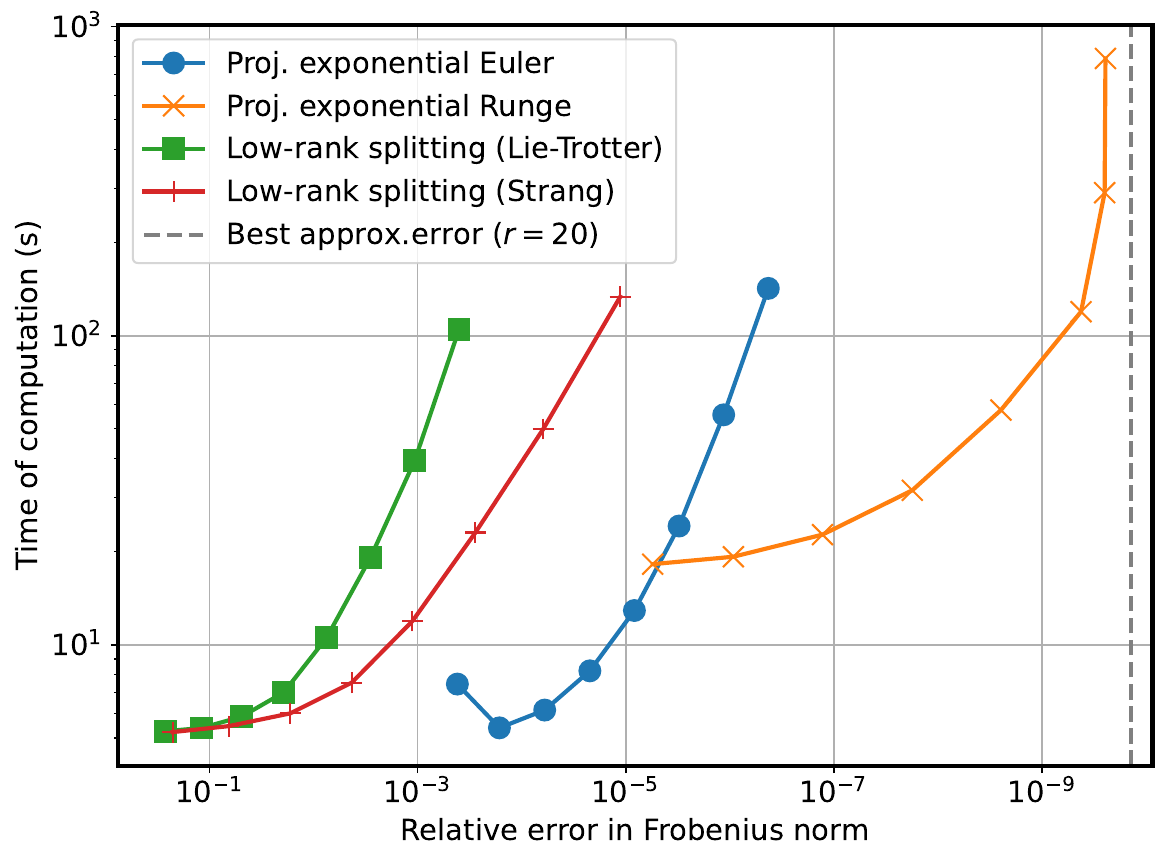}
        \subcaption{Performance. The dotted line is the approximation error that cannot be outperformed.}
    \end{subfigure}
    \caption{The new methods applied to the differential Riccati equation. The low-rank splitting methods are from \cite{ostermann2019convergence}.}
    \label{fig: riccati global error}
\end{figure}

Figure \ref{fig: riccati global error} shows a comparison between the new methods and the existing methods described in \cite{ostermann2019convergence}. 
The figure on the left shows the global error. As we can see, projected exponential Euler and Runge have convergence of order one and two, respectively. 
That is in agreement with Theorem \ref{theorem: projected exponential Euler convergence} and Theorem \ref{theorem: convergence of exponential Runge}.
Again, we observe that the new methods have smaller hidden constants compared to existing methods.
With enough time steps, the projected exponential Runge method reaches the approximation error which is close to $10^{-10}$.
On the other hand, the low-rank splitting (Strang) does not show a clear convergence of order two when the step size is too large.
We also provide the time of computation on the right figure.
For any given accuracy, we can see that the new methods are faster than existing methods.

\subsubsection*{Krylov approximation errors}

Let us consider the three different Krylov spaces described in Section \ref{sec: efficient implementation}.
Since the spectrum of $A$ is unbounded under mesh refinements, we use a single repeated pole $k / \sqrt{2},$ which turns out to be a simple but already good choice of poles, see \cite[Section 4.2]{guttel2013rational}.
Strictly speaking, the theoretical bounds derived in Section \ref{sec: efficient implementation} are only true when optimal poles are used. 
Nevertheless, the next figure includes the bounds for visual reference.

\begin{figure}[!ht]
    \begin{subfigure}{0.49\textwidth}
        \includegraphics[width=\textwidth]{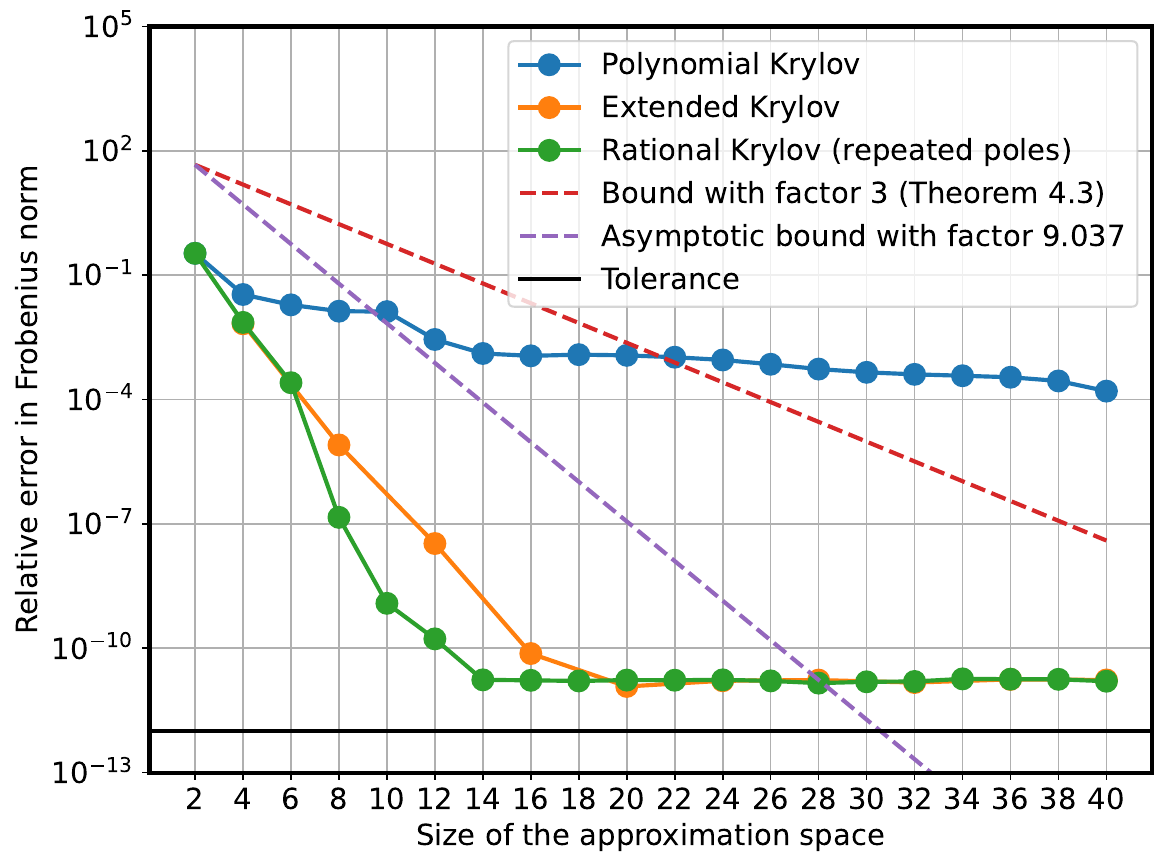}
        \subcaption{First order approximation: relative errors between the full system \eqref{eq: equivalence for projected Euler} and the reduced system \eqref{eq: reduced Sylvester IVP}.}
    \end{subfigure}
    \begin{subfigure}{0.49\textwidth}
        \includegraphics[width=\textwidth]{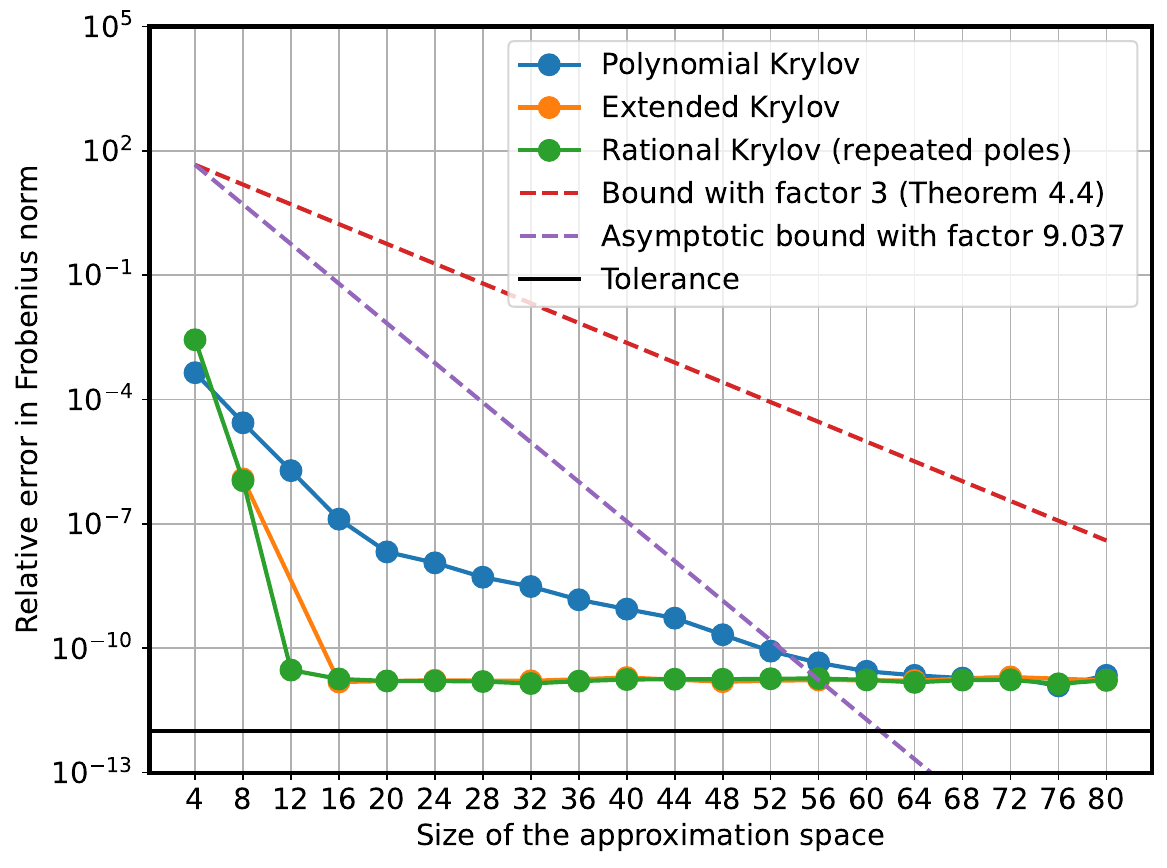}
        \subcaption{Second order approximation: relative errors between the full system \eqref{eq: system order 2} and the reduced system \eqref{eq: reduced system order 2}.}
    \end{subfigure}
    \caption{Approximation error made by the Krylov techniques for reducing an ODE.}
    \label{fig: krylov approximation error}
\end{figure}

Figure \ref{fig: krylov approximation error} shows the error made by each type of Krylov space as the size of the space grows, at the time $t=0.01$. 
In total, we perform $k=20$ iterations.
We start from the same initial value $X_0$ truncated to rank~$1$. 
Therefore, the approximation space starts from size~$2$ for the first order approximation, and size~$4$ for the second order approximation. 
We deliberately choose a very small rank since it produces clearer convergence plots. 
In our experiments, a higher rank of approximation always resulted in faster convergence than theoretically predicted.
As we can see, in both cases, the rational Krylov shows the faster convergence per iteration.
At the opposite, the polynomial Krylov space converges slowly, which implies large dimensional approximation space, and therefore large memory footprint. 
The extended Krylov space is somewhere in between and seems to be a good practical choice since it avoids the problem of choosing the poles that can be complex numbers (as in the case of optimal poles).

\subsection{Application to the Allen-Cahn equation}

The following example is inspired from \cite{rodgers2023implicit} and describes the process of phase separation in multi-component alloy systems \cite{allen1972ground, allen1973correction}. In its simplest form, the partial differential equation is given by
\begin{equation}
\frac{\partial f}{\partial t} = \varepsilon \Delta f + f - f^3,
\end{equation}
where $\Delta f$ is the diffusion term and $f - f^3$ is the reaction term. 
The stiffness is controlled with the small parameter $\varepsilon$.
Similarly to the paper cited above, we consider the initial condition
\begin{equation*}
f_0(x,y) = \frac{\left[ e^{-\tan^2(x)} + e^{-\tan^2(y)} \right] \sin(x) \sin(y)}{1 + e^{| \mathrm{csc}(-x/2) |} + e^{| \mathrm{csc}(-y/2)|} }. 
\end{equation*}
The space is $\Omega = [0, 2 \pi]^2$ and the mesh contains $256 \times 256$ points. 
The time interval is $[0, 10]$.
The discretization is done via finite differences, leading to a matrix differential equation of the form
\begin{equation} \label{eq: Allen-Cahn}
\dt{X}(t) = A X(t) + X(t) A + X(t) - X(t)^3, \quad (X_0)_{ij} = f_0(x_i, y_j),
\end{equation}
where the power is to be taken element-wise ($X^3 = X * X * X$ where $*$ is the Hadamard product).
The numerical rank of the initial value is $28$.
Because of the Hadamard product, the non-linearity is clearly more challenging than in the Lyapunov and Riccati differential equations.
A method for performing efficient Hadamard product with low-rank matrices is described in \cite{kressner2017recompression}, which avoids forming the dense matrix.
The growth in rank is therefore limited since $\rank(X^3) \leq \rank(X)^3$.
The computations are thus feasible when the rank $r$ is (very) small compared to the size $n$ of the matrix. 
In this example, we take $r=2$.

\begin{figure}[!ht]
\centering
    \begin{subfigure}{0.25\textwidth}
        \includegraphics[width=0.8\textwidth]{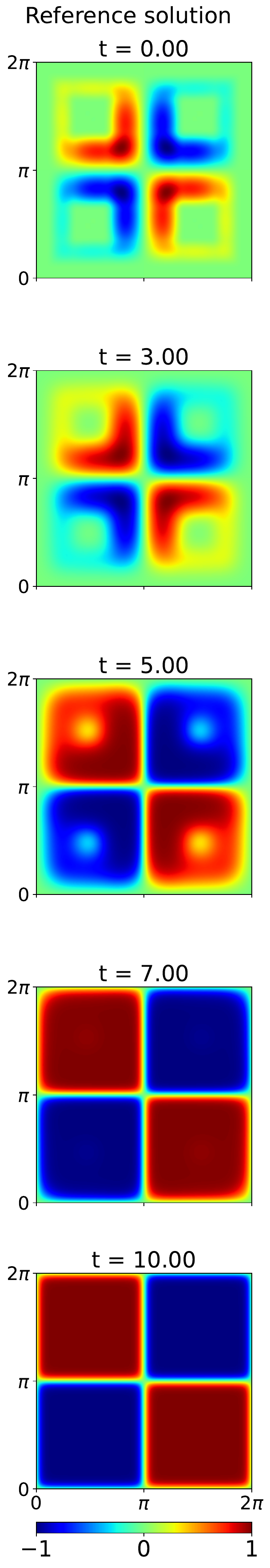}
        \subcaption{Comput. time: $\approx 5.26s$.}
    \end{subfigure}
    \begin{subfigure}{0.25\textwidth}
        \includegraphics[width=0.8\textwidth]{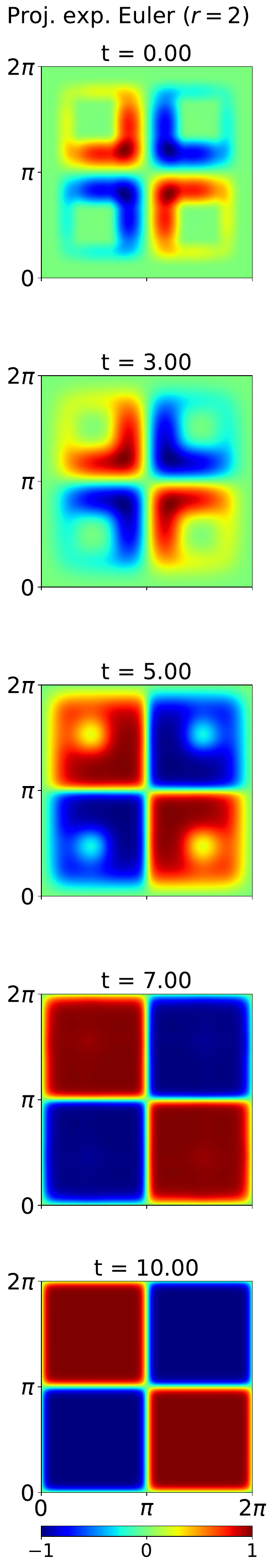}
        \subcaption{Comput. time: $\approx 0.22s$.}
    \end{subfigure}
    \begin{subfigure}{0.25\textwidth}
        \includegraphics[width=0.8\textwidth]{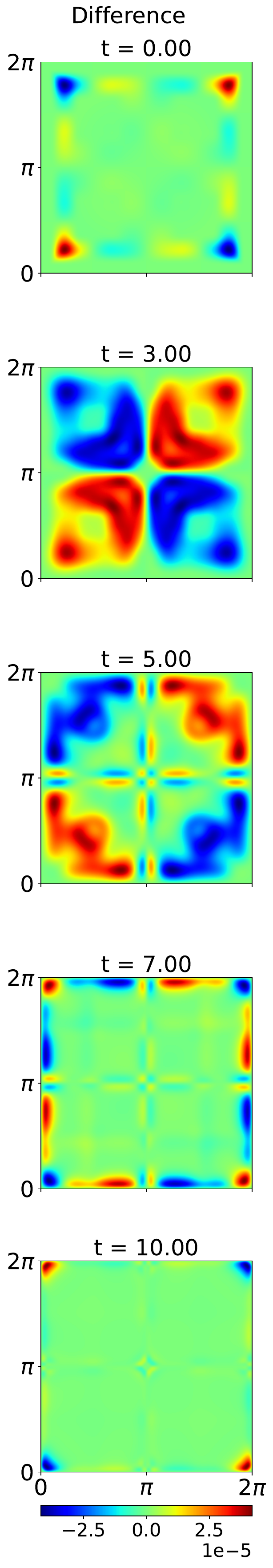}
        \subcaption{Normalized error.}
    \end{subfigure}
    \caption{Solution to Allen-Cahn equation \eqref{eq: Allen-Cahn} at several time points. The reference solution is computed with scipy's default method RK45. The computation time is the total time of computation for all steps. The projected exponential Euler method give excellent qualitative prediction and is significantly cheaper.}
    \label{fig: Allen-Cahn over time}
\end{figure}

In Figure \ref{fig: Allen-Cahn over time}, we give a qualitative comparison between the reference solution and the approximation made by projected exponential Euler. 
Our method already gives an excellent qualitative approximation despite the small rank $r=2$.
Moreover, the method is significantly faster than the reference solver, which here is scipy's solver RK45 with default tolerance ($10^{-8}$).

\subsection{Rank-adaptive methods}

Often in practice, the rank of approximation is not known a priori. 
In such a situation, we would like to give a tolerance and let the algorithm choose automatically the appropriate rank at each time step. 
Such methods are called rank-adaptive, and several of them have already been proposed, see \cite{ceruti2022rank}, \cite{hochbruck2023rank}.

Fortunately, the techniques proposed in this paper are easy to transform into rank-adaptive techniques.
Indeed, each truncation can be performed at any given tolerance instead of a fixed rank. 
For example, the adaptive version of projected exponential Euler is given by the iteration
\begin{align*} \label{eq: rank-adaptive projected exponential Euler}
    Y_1^{\mathrm{APE}} = \mathcal{T}_{\tau} \left( e^{h \cL} Y_0 + h \varphi_1 (h \cL) \proj{Y_0}{\cG (Y_0)} \right),
\end{align*}
where $\mathcal{T}_{\tau}$ is the truncated SVD up to tolerance $\tau$.
It naturally extends to any higher-order scheme.

\begin{figure}[!ht]
    \centering
    \begin{subfigure}{0.32\textwidth}
    \includegraphics*[width=\textwidth]{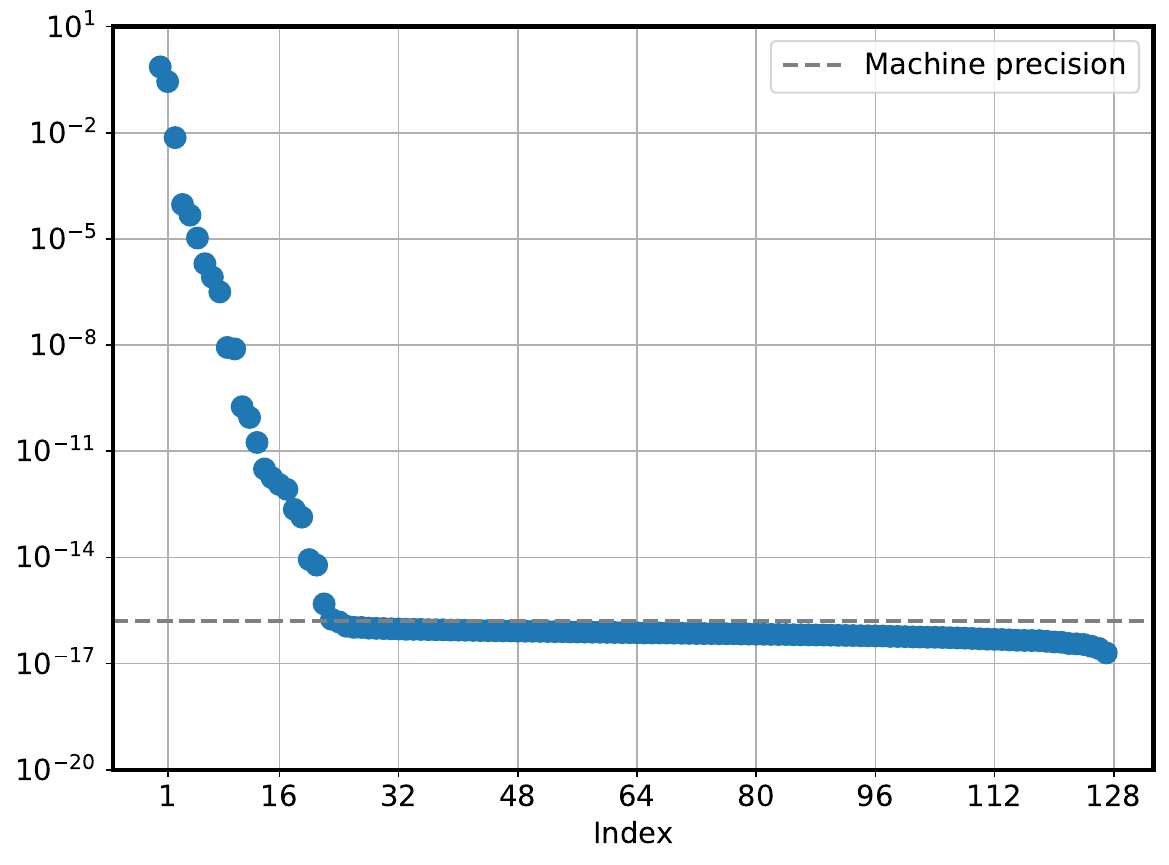}
    \end{subfigure}
    \begin{subfigure}{0.32\textwidth}
    \includegraphics*[width=\textwidth]{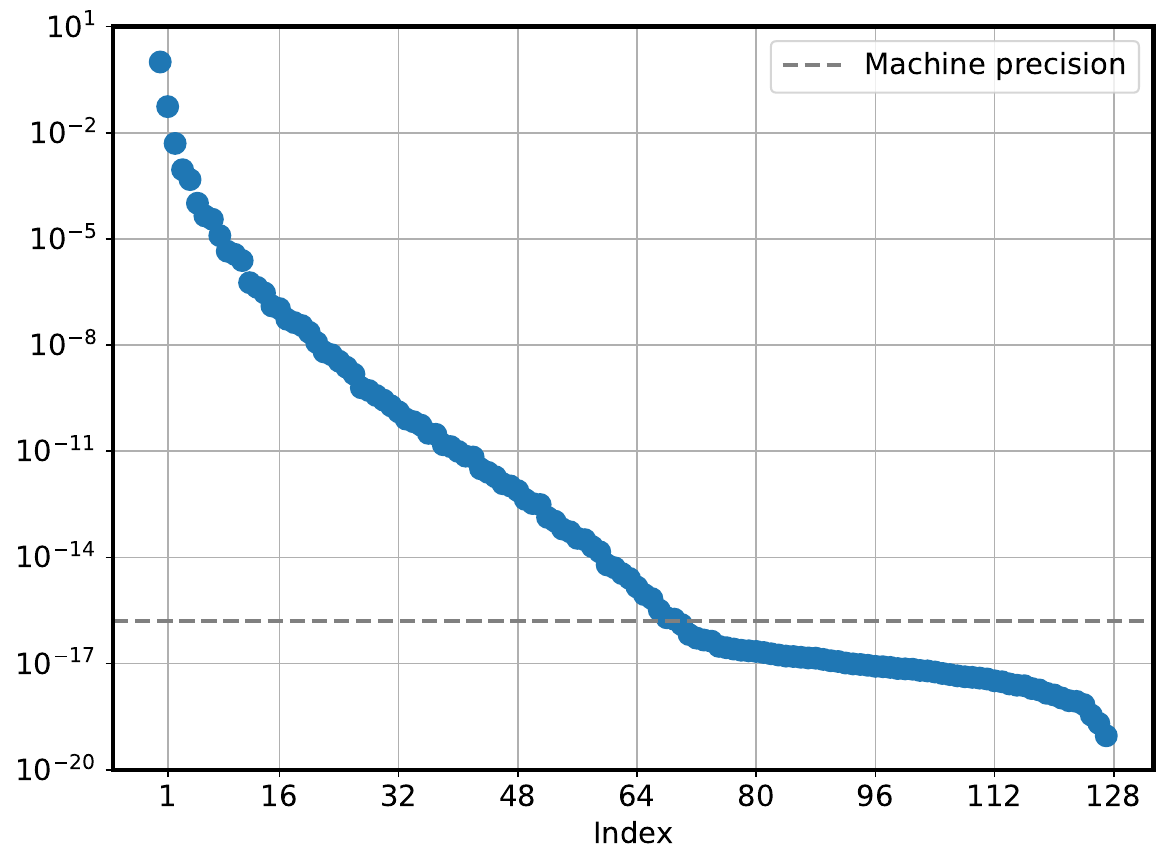}
    \end{subfigure}
    \begin{subfigure}{0.32\textwidth}
    \includegraphics*[width=\textwidth]{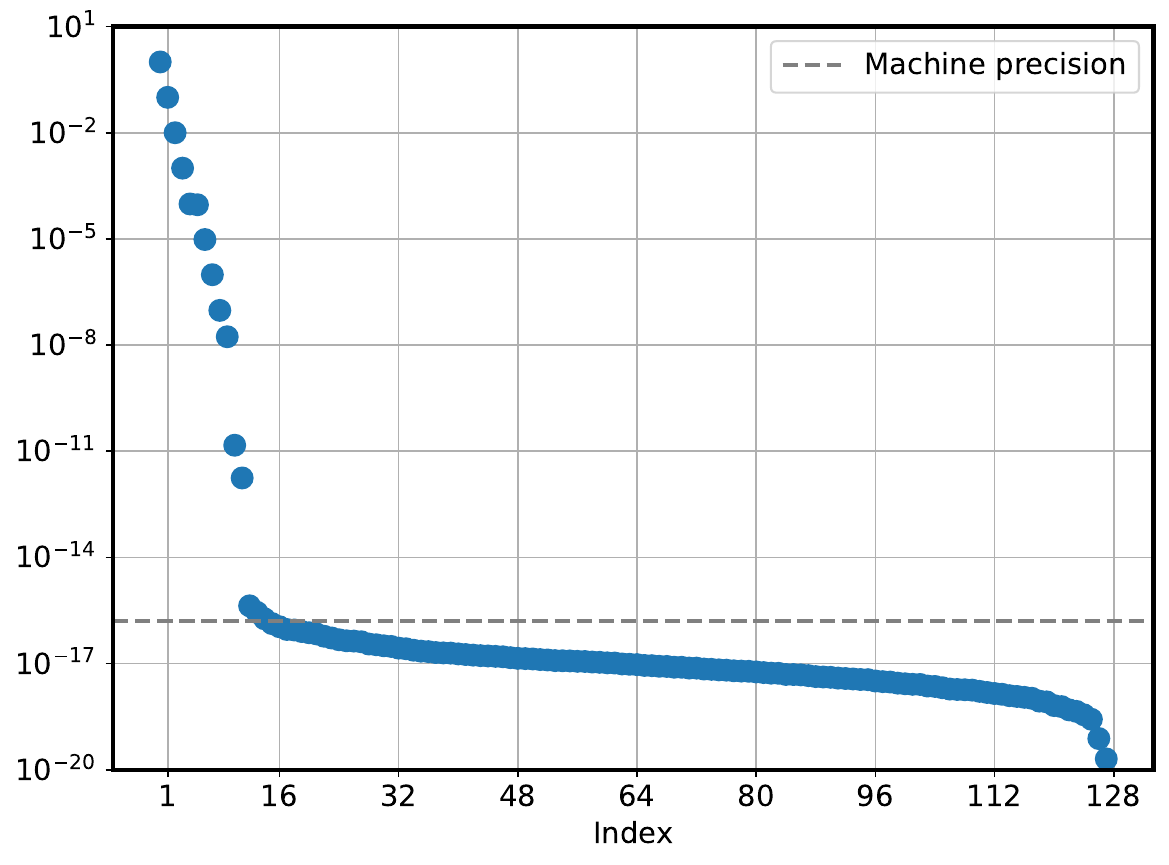}
    \end{subfigure}
    \begin{subfigure}{0.32\textwidth}
    \includegraphics*[width=\textwidth]{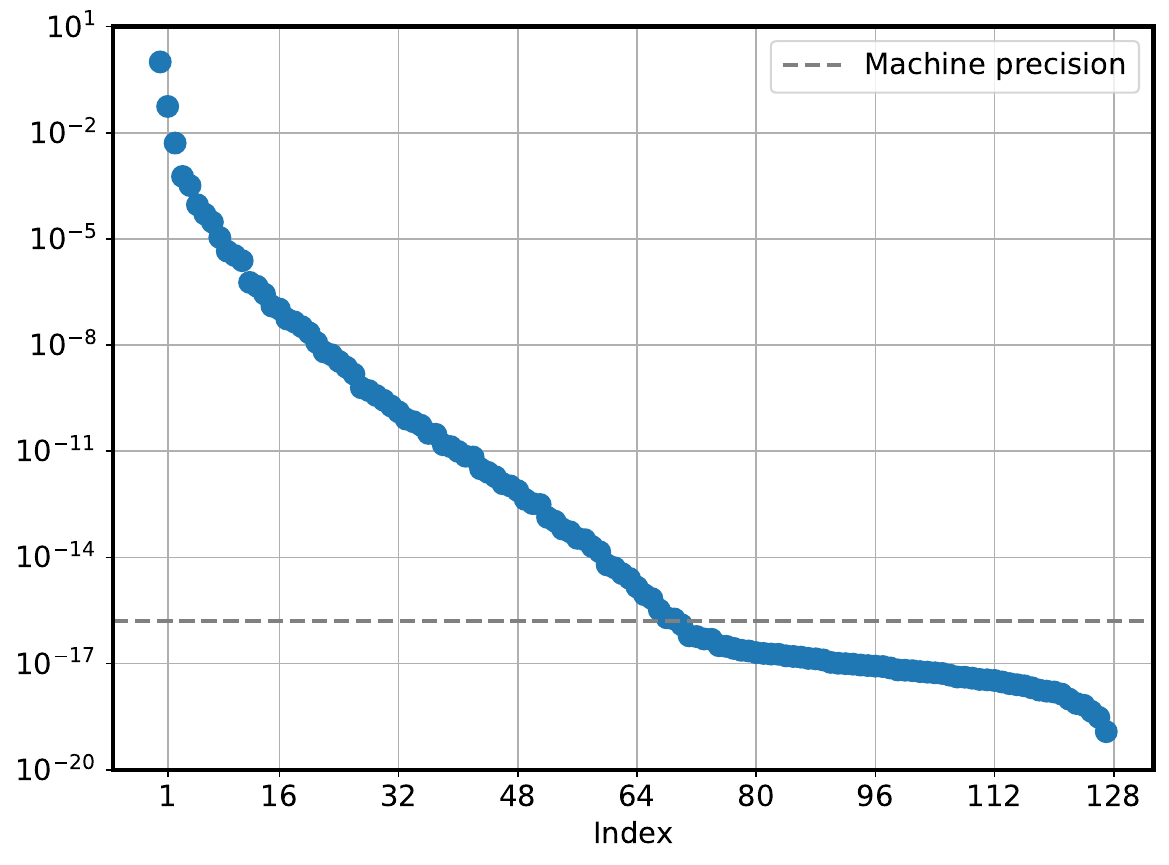}
    \end{subfigure}
    \begin{subfigure}{0.32\textwidth}
    \includegraphics*[width=\textwidth]{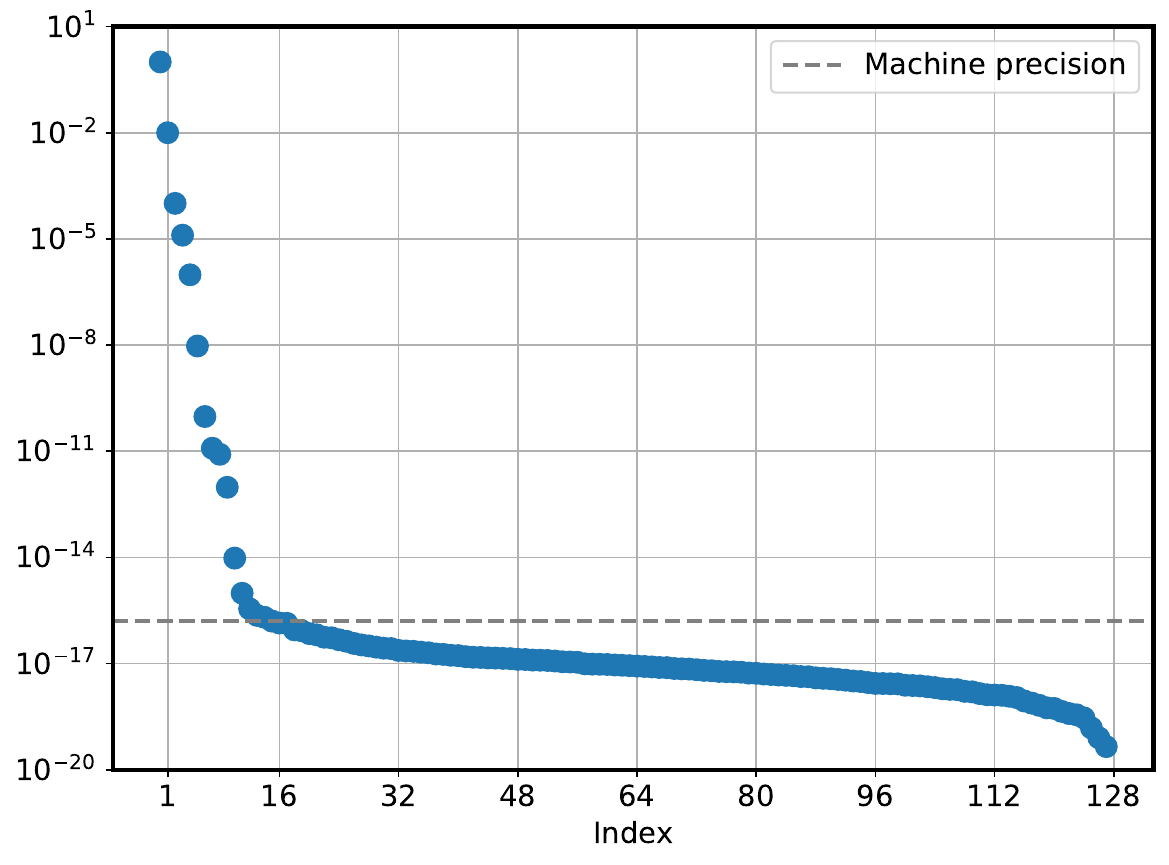}
    \end{subfigure}
    \caption{Singular values of the reference solution at several time points.}
    \label{fig: singular values over time for adaptive experiments}
\end{figure}

Let us consider again the Lyapunov differential equation \eqref{eq:lyapunov_diffeq} with a time-dependent source. 
In order to verify the capabilities of the method, we choose a time-dependent source that will modify the singular values significantly so that the ideal rank for a given tolerance is often changing:
\begin{equation*}
C(t) = \left\{ 
\begin{aligned}
&A X_1 + X_1 A^T, \quad &&\text{if } 0 \leq t < 0.2,  \\
&A X_1 + X_1 A^T + (t - 0.2)/0.2 \cdot (A X_2 + X_2 A^T - A X_1 + X_1 A^T), \quad &&\text{if } 0.2 \leq t < 0.4, \\
&A X_2 + X_2 A^T, \quad &&\text{if } 0.4 \leq t < 0.6, \\
&A X_2 + X_2 A^T + (t - 0.6)/0.2 \cdot (A X_1 + X_1 A^T - A X_2 + X_2 A^T), \quad &&\text{if } 0.6 \leq t < 0.8, \\
&A X_1 + X_1 A^T, \quad &&\text{if } 0.8 \leq t \leq 1,
\end{aligned}
\right.
\end{equation*}
where $X_1 = Q \Sigma_1 Q^T$ and $X_2 = Q \Sigma_2 Q^T$ with singular values $\Sigma_1=\mathrm{diag}(1, 10^{-2}, 10^{-4}, \ldots, 10^{-12}, 10^{-14})$, $\Sigma_2 = \mathrm{diag}(1, 10^{-1}, 10^{-2}, \ldots, 10^{-7}, 10^{-8})$ and random $Q \in \R^{n \times 9}$ with othonormal columns.
As we can see in Figure \ref{fig: singular values over time for adaptive experiments}, the singular values of the reference solution have essentially five phases, one for each time interval.

\begin{figure}[!ht]
    \centering
    \begin{subfigure}{0.49\textwidth}
        \includegraphics*[width=\textwidth]{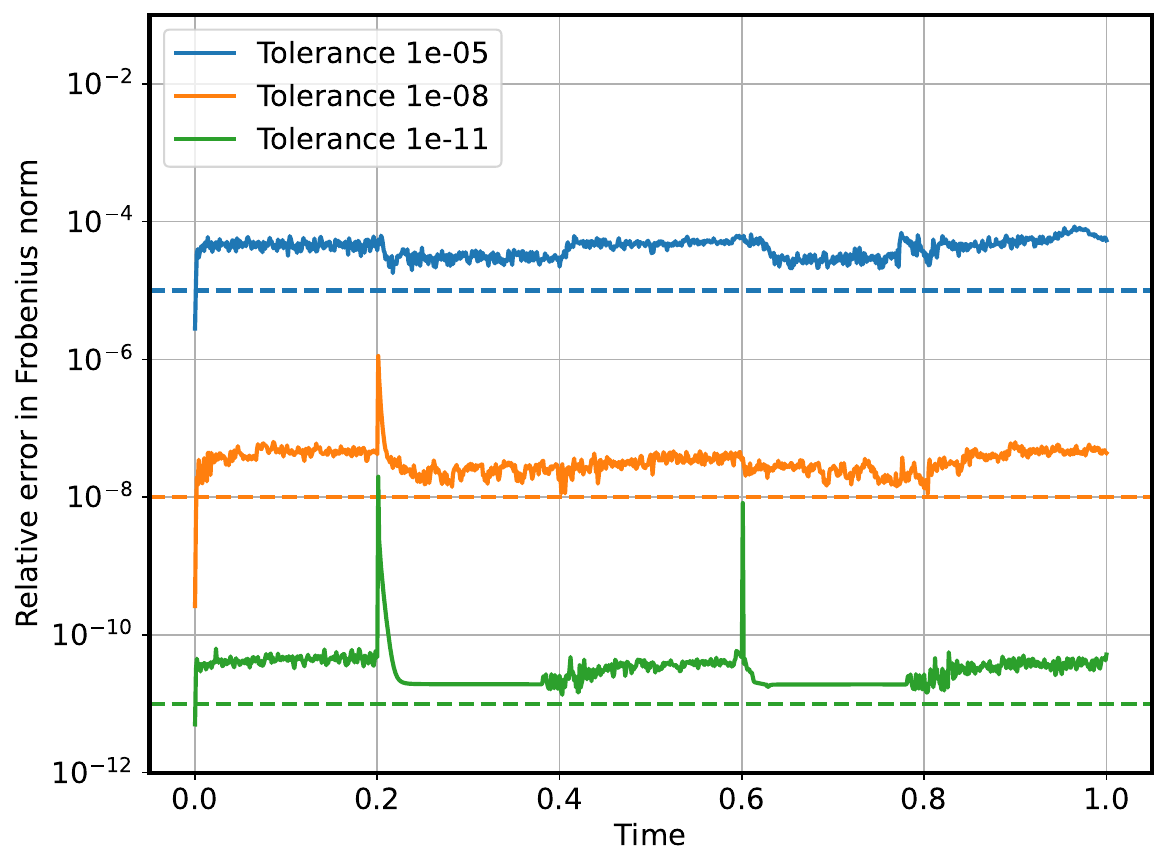}
        \subcaption{Relative error over time.}
    \end{subfigure}
    \begin{subfigure}{0.49\textwidth}
        \includegraphics[width=\textwidth]{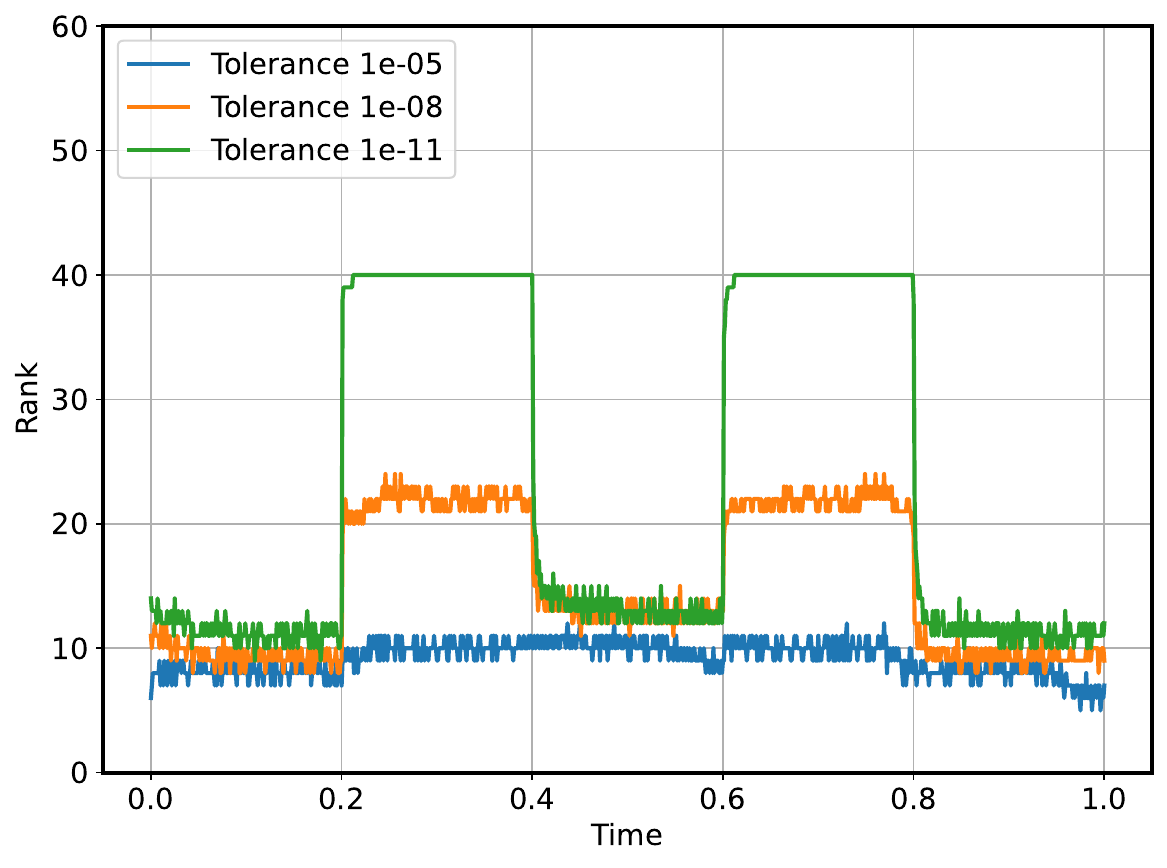}
        \subcaption{Rank evolution over time.}
    \end{subfigure}
    \caption{Adaptive projected exponential Runge applied to the Lyapunov differential equation with time-dependent source. Step size $h = 0.001$.}
    \label{fig: adaptive lyapunov}
\end{figure}

In Figure \ref{fig: adaptive lyapunov}, we applied the rank-adaptive projected exponential Runge method with three different tolerances. 
The step size is $h=0.001$, and we again used extended Krylov spaces with one iteration. 
As we can see, the method is able to rapidly adapt the rank, and the small tolerances are able to capture the five phases of the singular values. 
The peaks in the errors are most likely due to the discontinuous nature of the problem. 
The relative errors are close to the given tolerances as it should be. 
Overall, the adaptive methods show good results with no extra cost due to the adaptive rank.

\section{Conclusion and future works} \label{sec:conclusion}

In a nutshell, we proposed a new class of methods for the dynamical low-rank approximation that we called projected exponential methods. 
In particular, we focused on two methods, namely projected exponential Euler and Runge, which we have proven to be robust-to-stiffness with a stiff order one and two, respectively, up to the modelling error due to low-rank. 
In addition, we described a third, intermediate, technique that only requires evaluating $\varphi_1$ functions and that has good practical convergence properties.
In order to apply the schemes, we proposed Krylov techniques for efficient implementation, reducing the memory footprint and accelerating the computations.
Those techniques are also theoretically analyzed, leading to convergence bounds that can be a guide when choosing the size of the approximation space.
Finally, numerical experiments verify the theoretical derivations and show the good performance of the new methods.

An obvious future work is to extend the analysis to higher-order methods. 
However, higher-order methods will require a significantly larger approximation space, making the computations not feasible in most cases.
Recently, many efforts are made into developing iterative methods that are low-rank compatible.
It is possible that higher-order projected exponential methods become viable in the near future.

Another possible work would be to extend the framework to dynamical tensor approximation \cite{koch2010dynamical}. Certain discretized high-dimensional Schrödinger equations that are stiff are numerically challenging to solve in a large-scale setting. Interestingly, the equation has the structure required for projected exponential integrators, and is therefore a good motivation for studying further those methods in high dimensions.

\backmatter
\bibliography{perk_arxiv_v2}

\appendix

\section{Convergence proof of exponential Euler} \label{appendix: exponential euler}

The following is the proof of Theorem \ref{theorem: convergence of exponential Euler}.

\begin{proof}
Let us start with a Taylor expansion of the closed form solution:
\begin{align*}
X(t_n) &= e^{h \cL} X(t_n) + \int_0^h e^{(h-\tau) \cL} f(t_n + \tau) d\tau \\
&= e^{h \cL} X(t_n) + h \varphi_1(h \cL) f(t_n) + \int_0^h e^{(h-\tau) \cL} \int_0^{\tau} f'(t_n + \sigma) d \sigma d\tau.
\end{align*}
Therefore, the error $E_{n+1} = X_{n+1}^{\mathrm{E}} - X(t_{n+1})$ verifies the following recursion:
$$E_{n+1} = e^{h \cL} E_n + h \varphi_1(h \cL) \left[\cG(X_n) - f(t_n) \right] - \delta_{n+1},$$
where $\delta_{n+1}$ is the remainder term at the step $n+1$.
Taking the norm and using the Lipschitz continuity of $\cG$ gives that
$$\norm{E_{n+1}} \leq e^{h \ell} \norm{E_n} + h \cdot \varphi_1(h \ell) L_{\cG} \cdot \norm{E_n} + \norm{\delta_{n+1}}$$
which leads to
$$\norm{E_{n+1}} \leq (1 + h L_*) \norm{E_n} + \norm{\delta_{n+1}},$$
where $L_* = \ell + \varphi_1(h \ell) L_{\cG}$.
The remainder term is bounded by
\begin{align*}
\norm{\delta_{n+1}} &\leq \int_0^h e^{(h-\tau) \ell} \tau d\tau \max_{0 \leq s \leq h} \norm{f'(t_n + s)} \\
&\leq h^2 \varphi_2(h \ell) \max_{0 \leq s \leq h} \norm{f'(t_n+s)} \\
&\leq h^2 \varphi_2(h \ell) \max_{0 \leq s \leq t_{n+1}} \norm{f'(s)} = h \beta.
\end{align*}
Solving the recursion as done in \eqref{eq: solving recursion} leads to
\begin{align*}
\norm{X_n - X(t_n)} &\leq \frac{e^{nhL_*} - 1}{L_*} \cdot \varphi_2(h \ell) \cdot \max_{0 \leq t \leq t_n} \norm{f'(t)}.
\end{align*}
For  all $h \leq h_0$ small enough, the factors above reduce to a constant $C(\ell, L_{\cG}, h_0)$ as stated in the theorem.
\end{proof}

\end{document}